\documentclass[12pt,reqno]{amsart}
\usepackage{amssymb,mathrsfs,eqnarray,bbold,amsxtra,xcolor,hyperref} 

\theoremstyle{plain}
\newtheorem{theorem}{Theorem}[section]
\newtheorem{proposition}[theorem]{Proposition}
\newtheorem{corollary}[theorem]{Corollary}
\newtheorem{lemma}[theorem]{Lemma}

\theoremstyle{definition}

\newtheorem{example}[theorem]{Example}

\numberwithin{equation}{section}

\DeclareMathOperator{\rad}{Rad}
\DeclareMathOperator{\dist}{dist}

\DeclareMathOperator{\tr}{tr}

\DeclareMathOperator{\card}{card}
\renewcommand{\le}{\leqslant}
\renewcommand{\ge}{\geqslant}

\begin{document}
\title[On operators with infinitely many spectral components]{On operators with infinitely many spectral components}
\author[J. Bra\v{c}i\v{c}]{Janko Bra\v{c}i\v{c}}
\address{Faculty of Natural Sciences and Engineering, University of Ljubljana,
A\v{s}ker\v{c}eva c. 12, SI-1000 Ljubljana, Slovenia}
\email{janko.bracic@ntf.uni-lj.si}
\author[M. Kandi\'c]{Marko Kandi\'c}
\address{Faculty of Mathematics and Physics, University of Ljubljana,
Jadranska 19, 1000 Ljubljana, Slovenia \ \ \ and \ \ \ Institute of Mathematics, Physics, and Mechanics,
Jadranska 19, 1000 Ljubljana, Slovenia}
\email{marko.kandic@fmf.uni-lj.si}
\begin{abstract}
For a unital complex Banach algebra $\mathscr A$ and $n\in\mathbb N$, let $\mathcal C_n(\mathscr A)$ be the set of all elements 
whose spectrum has at least $n$ connected components, and let 
$\mathcal C_\infty(\mathscr A)=\bigcap_{n\in\mathbb N}\mathcal C_n(\mathscr A)$. We prove that 
$\mathcal C_\infty(\mathscr A)\neq\varnothing$ if and only if $\mathscr A$ contains a sequence of non-zero pairwise orthogonal 
idempotents, and that the number of components of a spectrum is always finite, $\aleph_0$ or $\mathfrak c$. This allows us to extend 
the spectral depth $d_\sigma(\mathscr A)$ to a cardinal-valued invariant with values in $\mathbb N\cup\{\aleph_0,\mathfrak c\}$. We 
show that $d_\sigma(\mathscr A)=\mathfrak c$ precisely when $\mathscr A$ contains a dyadic family of idempotents indexed by the 
full binary tree. In that case, some element has spectrum homeomorphic to the Cantor space $\mathbf 2^{\mathbb N}$. The 
supremum defining $d_\sigma(\mathscr A)$ is attained when the spectral classes of $\mathscr A$ are dense and when $\mathscr A$ is 
commutative, but not in general. We construct a unital complex Banach algebra with $d_\sigma(\mathscr A)=\aleph_0$ and 
$\mathcal C_\infty(\mathscr A)=\varnothing$. If $\Delta$ is the character space of a commutative $\mathscr A$ and $\pi_0(\Delta)$ is 
the space of all components of $\Delta$, then $d_\sigma(\mathscr A)=\card\bigl(\pi_0(\Delta)\bigr)$ when $\pi_0(\Delta)$ is finite, 
$d_\sigma(\mathscr A)=\aleph_0$ when $\pi_0(\Delta)$ is infinite and scattered, and $d_\sigma(\mathscr A)=\mathfrak c$ when 
$\pi_0(\Delta)$ is not scattered. We also prove that the spectral depth is additive with respect to finite direct sums and that the 
spectral classes are invariant under injective holomorphic functions.
\end{abstract}
\maketitle
\section{Introduction} \label{sec:introduction}
The spectrum of an element of a complex Banach algebra is a non-empty compact subset of the complex plane, and the number of its 
connected components is the crudest, and yet most robust, invariant attached to it. In the companion paper \cite{BK26} we organized a 
complex Banach algebra $\mathscr A$ with unit $e$ into the decreasing chain of \emph{spectral classes}
$$\mathcal C_n(\mathscr A)=\{a\in\mathscr A;\;\sigma(a) \text{ has at least } n \text{ components}\} \qquad (n\in\mathbb N),$$
and we studied the \emph{spectral depth}
$$d_\sigma(\mathscr A)=\sup\{n\in\mathbb N;\;\mathcal C_n(\mathscr A)\neq\varnothing\}\in\mathbb N\cup\{\infty\}.$$
The main structural result there, \cite[Proposition~3.1]{BK26}, is that the whole hierarchy is governed by idempotents: 
$\mathcal C_n(\mathscr A)\neq\varnothing$ if and only if the unit of $\mathscr A$ decomposes into $n$ non-zero pairwise orthogonal 
idempotents.

The present paper is devoted to the limiting class
$$\mathcal C_\infty(\mathscr A)=\bigcap_{n=1}^{\infty}\mathcal C_n(\mathscr A)$$
of elements whose spectrum has infinitely many components, and to the finer classification of such elements according to how many 
components their spectra have.

Section~\ref{sec:Cinfty-spectral-depth} contains the counterpart, at the level $n=\infty$, of
\cite[Proposition~3.1]{BK26}. Proposition~\ref{prop:Cinfty-idempotent-characterization} states that 
$\mathcal C_\infty(\mathscr A)\neq\varnothing$ if and only if $\mathscr A$ contains a sequence of non-zero pairwise orthogonal 
idempotents, if and only if some element of $\mathscr A$ has spectrum of the form $\{0\}\cup\{\lambda_j;\;j\in\mathbb N\}$, where 
$(\lambda_j)_{j\in\mathbb N}$ is a sequence of distinct non-zero complex numbers converging to $0$. The passage from finite to 
infinite is, however, not automatic. The condition $d_\sigma(\mathscr A)=\infty$ says that every class $\mathcal C_n(\mathscr A)$ is
non-empty, whereas $\mathcal C_\infty(\mathscr A)\neq\varnothing$ says that their intersection is non-empty, and a decreasing chain 
of non-empty sets may have empty intersection. By Proposition~\ref{prop:Cinfty-idempotent-characterization}, asking whether the first 
condition implies the second amounts to the purely algebraic question whether an algebra containing, for every $n$, a family of $n$ 
non-zero pairwise orthogonal idempotents contains an infinite such family. Proposition~\ref{prop:dense-classes} gives an
affirmative answer for every algebra all of whose spectral classes are dense. By \cite[Theorem~4.7]{BK26} this covers $\mathscr B(X)$ 
for an infinite-dimensional complex Banach space $X$.

Section~\ref{sec:orthogonal-infinitude} shows that the answer is negative in general. We construct a unital complex Banach algebra of 
continuous $2\times2$ matrix-valued functions, subject to a boundary condition built from a coordinate-swap map on a product of 
finite sets, whose identity decomposes into $n$ non-zero pairwise orthogonal idempotents for every $n$, while the algebra contains no 
sequence of non-zero pairwise orthogonal idempotents (Corollary~\ref{cor:infinite-spectral-depth-not-attained}). Thus 
$\mathcal C_n(\mathscr A)\neq\varnothing$ for every $n\in\mathbb N$, although $\mathcal C_\infty(\mathscr A)=\varnothing$.

Section~\ref{sec:cardinality-spectral-components} determines the possible cardinalities. Since the set $\pi_0(\Omega)$ of all 
components of a compact metrizable space $\Omega$ is again compact, metrizable and totally disconnected, hence Polish, when 
equipped with the quotient topology, the Cantor--Bendixson theorem yields Theorem~\ref{thm:possibility-for-cardinality}: the number 
of components of a spectrum belongs to $\mathbb N\cup\{\aleph_0,\mathfrak c\}$, with no intermediate cardinality possible. This is a 
theorem of ZFC, independent of the Continuum Hypothesis. Accordingly, $\mathcal C_\infty(\mathscr A)$ splits into the two disjoint 
classes $\mathcal C_{\aleph_0}(\mathscr A)$ and $\mathcal C_{\mathfrak c}(\mathscr A)$, and the spectral depth extends to a 
cardinal-valued invariant $d_\sigma(\mathscr A)\in\mathbb N\cup\{\aleph_0,\mathfrak c\}$. The main result of the section,
Theorem~\ref{thm:dyadic-idempotents-Cantor-spectrum}, is the exact analogue of \cite[Theorem~3.4]{BK26} one level higher: if 
$\mathscr A$ contains a family $\{u_s\}$ of non-zero idempotents indexed by the full binary tree and satisfying $u_\varnothing=e$, 
$u_s=u_{s0}+u_{s1}$ and $u_{s0}u_{s1}=0$, then $\mathscr A$ contains an element whose spectrum is homeomorphic to the Cantor 
space $\mathbf 2^{\mathbb N}$. Together with a Cantor--Bendixson argument, this gives 
Corollary~\ref{cor:continuum-spectral-class-tree-characterization}, a characterization of the non-emptiness of 
$\mathcal C_{\mathfrak c}(\mathscr A)$ by such dyadic families, and Corollary~\ref{cor:depth-trichotomy}, a trichotomy for the spectral 
depth in terms of the idempotent structure of $\mathscr A$ alone. Both values $\aleph_0$ and $\mathfrak c$ really occur for algebras 
of operators. For an infinite-dimensional complex Banach space $X$, we have $d_\sigma(\mathscr B(X))=\mathfrak c$ whenever $X$ 
contains a complemented subspace with an unconditional Schauder basis (Proposition~\ref{prop:unconditional-basis}), whereas
$d_\sigma\bigl(\mathscr B(X_{\mathrm{AH}})\bigr)=\aleph_0$ for the Argyros--Haydon space (Example~\ref{ex:AH}). 
Section~\ref{sec:cardinality-spectral-components} also contains two permanence properties. The classes $\mathcal C_\nu(\mathscr A)$ 
are invariant under injective holomorphic functions (Proposition~\ref{prop:functional-calculus}). The spectral depth is additive with 
respect to finite direct sums, in the sense of cardinal arithmetic (Theorem~\ref{thm:depth-additive}).

In Section~\ref{sec:commutative-banach-algebras} we show that $d_\sigma(\mathscr A)=\infty$ implies 
$\mathcal C_\infty(\mathscr A)\neq\varnothing$ for a commutative algebra $\mathscr A$. The \v{S}ilov idempotent theorem identifies 
the idempotents of $\mathscr A$ with the clopen subsets of $\Delta$, the character space of $\mathscr A$, hence with the clopen 
subsets of $\pi_0(\Delta)$. The results of Sections~\ref{sec:Cinfty-spectral-depth} and \ref{sec:cardinality-spectral-components} 
translate the three regimes into topological properties of $\pi_0(\Delta)$. Theorem~\ref{thm:commutative-depth} states that 
$d_\sigma(\mathscr A)=\card\bigl(\pi_0(\Delta)\bigr)$ if $\pi_0(\Delta)$ is finite, $d_\sigma(\mathscr A)=\aleph_0$ if $\pi_0(\Delta)$ is 
infinite and scattered, and $d_\sigma(\mathscr A)=\mathfrak c$ if $\pi_0(\Delta)$ is not scattered. In all three cases, the supremum is 
attained. In particular, $d_\sigma(\mathscr A)=d_\sigma\bigl(C(\Delta)\bigr)$, which means that the spectral depth of a commutative 
Banach algebra does not see its radical, and \cite[Proposition~4.5]{BK26} is recovered and completed. Moreover, in contrast to the 
general situation, a commutative algebra of infinite spectral depth is always orthogonally infinite. The section closes with the finite case, 
where finite spectral depth is shown to be equivalent to a decomposition of $\mathscr A$ into finitely many unital ideals without 
non-trivial idempotents (Theorem~\ref{thm:finite-depth-decomposition}), the decomposition being unique up to a permutation. 
The algebra $\mathbb M_n$ shows that commutativity cannot be dropped.

\section{Preliminaries} \label{sec:preliminaries}
Throughout the paper, $\mathscr A$ denotes a complex Banach algebra with unit $e$. We freely use the notation and the results of
\cite{BK26}. In this section we fix the terminology and collect the additional topological facts we shall need. By $\mathfrak c$ we
denote the cardinality of the continuum. For cardinals $\mu$ and $\nu$, the sum $\mu+\nu$ and the order relation $\mu\le \nu$ are
understood in the sense of cardinal arithmetic.

\subsection*{Components of a topological space}
Let $\Omega$ be a topological space. A subset of $\Omega$ is \emph{clopen} if it is both open and closed. A subset of $\Omega$ is
\emph{connected} if it cannot be written as the union of two non-empty disjoint relatively open subsets, and a \emph{connected
component}, or simply a \emph{component}, of $\Omega$ is a maximal connected subset of $\Omega$. The space $\Omega$ is
\emph{connected} if it has exactly one component, \emph{disconnected} otherwise, and \emph{totally disconnected} if every
component is a singleton. We denote by $\pi_0(\Omega)$ the set of all components of $\Omega$. Since every connected subset of
$\Omega$ meeting a clopen set is contained in it, each clopen subset of $\Omega$ is a union of components of $\Omega$.

A non-empty topological space is \emph{zero-dimensional} if it satisfies the separation axiom $T_1$ and has a base consisting of
clopen sets. By \cite[Theorem~6.2.1]{Eng89} every zero-dimensional space is totally disconnected and, by
\cite[Theorem~6.2.10]{Eng89}, every totally disconnected compact Hausdorff space is zero-dimensional. Thus, within the class of
compact Hausdorff spaces the two notions agree, and we use them interchangeably. We note that Engelking calls \emph{hereditarily
disconnected} what we call totally disconnected, and that in his book compact spaces are Hausdorff by definition. Recall also that in a compact Hausdorff space the
component of a point coincides with its quasi-component, that is, with the intersection of all clopen sets containing that point; see
\cite[Theorem~6.1.23]{Eng89}. We shall use repeatedly the following consequence, which is stated as Lemma~2.1 in \cite{BK26}. If
$n\in\mathbb N$ and a compact Hausdorff space $\Omega$ has at least $n$ components, then there are non-empty pairwise disjoint
clopen subsets $\Delta_1,\ldots,\Delta_n$ of $\Omega$ whose union is $\Omega$.

The following three lemmas are folklore. Since we could not find references, we include their proofs.

\begin{lemma}\label{lem:continuous-image-components}
Let $\phi\colon\Omega_1\to\Omega_2$ be a continuous mapping between topological spaces. Then
$\card\bigl(\pi_0\bigl(\phi(\Omega_1)\bigr)\bigr)\le\card\bigl(\pi_0(\Omega_1)\bigr)$.
\end{lemma}

\begin{proof}
Let $\Gamma\in\pi_0(\Omega_1)$. Since $\phi$ is continuous, $\phi(\Gamma)$ is a connected subset of $\phi(\Omega_1)$ and
therefore it is contained in a unique component $\Phi(\Gamma)$ of $\phi(\Omega_1)$. This defines a mapping
$\Phi\colon\pi_0(\Omega_1)\to\pi_0\bigl(\phi(\Omega_1)\bigr)$. To see that $\Phi$ is surjective, let
$\Delta\in\pi_0\bigl(\phi(\Omega_1)\bigr)$ be arbitrary. Pick $y\in \Delta$ and choose $x\in\Omega_1$ with $\phi(x)=y$. Let
$\Gamma$ be the component of $\Omega_1$ containing $x$. Then $\phi(\Gamma)$ is a connected set containing $y$. Since $\Delta$
is the component of $\phi(\Omega_1)$ containing $y$, we have $\phi(\Gamma)\subseteq\Delta$ and therefore
$\Phi(\Gamma)=\Delta$.
\end{proof}

\begin{lemma}\label{lem:infinitely-many-components}
Let $\Omega$ be a topological space with infinitely many components. Then there exists a sequence
$(\Delta_j)_{j\in\mathbb N}$ of non-empty pairwise disjoint clopen subsets of $\Omega$.
\end{lemma}

\begin{proof}
We construct the sets recursively. Put $\Omega_0=\Omega$. Suppose that $\Omega_{j-1}$ is a clopen subset of $\Omega$ with
infinitely many components. Then $\Omega_{j-1}$ is disconnected, which means that there is a non-empty proper clopen subset 
$\Gamma$ of $\Omega_{j-1}$. The sets $\Gamma$ and $\Omega_{j-1}\setminus \Gamma$ are non-empty disjoint clopen subsets of 
$\Omega_{j-1}$, and every component of $\Omega_{j-1}$ lies entirely in one of them. Hence at least one of them has infinitely many 
components. Denote such a set by $\Omega_j$ and the other one by $\Delta_j$. Since $\Omega_{j-1}$ is clopen in $\Omega$, both 
$\Omega_j$ and $\Delta_j$ are clopen in $\Omega$. By construction, $\Omega=\Omega_j\cup\Delta_1\cup\cdots\cup\Delta_j$ is a 
disjoint union. Each $\Delta_j$ is non-empty by the choice of $\Gamma$, and this decomposition shows that the $\Delta_j$ are 
pairwise disjoint.
\end{proof}

\begin{lemma}\label{lem:union-components}
Let $\Omega$ be a topological space and let $\Lambda_1$ and $\Lambda_2$ be non-empty closed subsets of $\Omega$. Then
\begin{equation}\label{eq:cardinality-components-inequality}
\card\bigl(\pi_0(\Lambda_1\cup\Lambda_2)\bigr)\le\card\bigl(\pi_0(\Lambda_1)\bigr)+\card\bigl(\pi_0(\Lambda_2)\bigr).
\end{equation}
If $\Lambda_1\cap\Lambda_2=\varnothing$, then
\begin{equation}\label{eq:cardinality-components-equality}
\card\bigl(\pi_0(\Lambda_1\cup\Lambda_2)\bigr)=\card\bigl(\pi_0(\Lambda_1)\bigr)+\card\bigl(\pi_0(\Lambda_2)\bigr).
\end{equation}
\end{lemma}

\begin{proof}
Let $\Gamma\in\pi_0(\Lambda_1)\cup\pi_0(\Lambda_2)$. Then $\Gamma$ is a connected subset of $\Lambda_1\cup\Lambda_2$,
and therefore it is contained in a unique component $\Phi(\Gamma)$ of $\Lambda_1\cup\Lambda_2$. This defines a mapping
$\Phi\colon\pi_0(\Lambda_1)\cup\pi_0(\Lambda_2)\to\pi_0(\Lambda_1\cup\Lambda_2)$. To see that $\Phi$ is surjective, let
$\Delta$ be a component of $\Lambda_1\cup\Lambda_2$ and let $x\in\Delta$. Then $x$ lies in $\Lambda_1$ or in $\Lambda_2$.
Let $\Gamma$ be the component of this set which contains $x$. Since $\Delta$ is the component of $\Lambda_1\cup\Lambda_2$
containing $x$, we have $\Gamma\subseteq\Delta$ and therefore $\Phi(\Gamma)=\Delta$. This proves
\eqref{eq:cardinality-components-inequality}.

Assume now that $\Lambda_1\cap\Lambda_2=\varnothing$. Each of the sets $\Lambda_1$ and $\Lambda_2$ is closed in
$\Lambda_1\cup\Lambda_2$ and it is the complement of the other one. Hence both are clopen in
$\Lambda_1\cup\Lambda_2$, and no connected subset of $\Lambda_1\cup\Lambda_2$ meets both of them. Let $\Gamma$ be a
component of $\Lambda_1$. Then $\Phi(\Gamma)$ is a connected set which meets $\Lambda_1$. Hence
$\Phi(\Gamma)\subseteq\Lambda_1$, and the maximality of $\Gamma$ gives $\Phi(\Gamma)=\Gamma$. The same holds for the
components of $\Lambda_2$. It follows that $\Phi$ is the identity on $\pi_0(\Lambda_1)\cup\pi_0(\Lambda_2)$, and therefore
$\pi_0(\Lambda_1\cup\Lambda_2)$ is the disjoint union of $\pi_0(\Lambda_1)$ and $\pi_0(\Lambda_2)$. This gives
\eqref{eq:cardinality-components-equality}.
\end{proof}

\subsection*{The topology of the space of components}
Let $\Omega$ be a compact Hausdorff space. The components of $\Omega$ partition $\Omega$, which means that the relation 
$x\sim y$, defined by $x, y\in\Omega$ lying in the same component, is an equivalence relation on $\Omega$. We equip the quotient 
$\pi_0(\Omega)$ with the quotient topology. As recalled above, components coincide with quasi-components, hence clopen sets 
separate distinct components. This is what lies behind the following theorem, which is \cite[Theorem~6.2.24]{Eng89} together with
\cite[Theorem~4.2.13]{Eng89}. We note that Engelking calls an equivalence relation \emph{closed} when the associated quotient map is closed.

\begin{theorem}\label{thm:component-space}
Let $\Omega$ be a compact Hausdorff space and let $\pi\colon \Omega\to\pi_0(\Omega)$ be the quotient map.
\begin{enumerate}
\item The space $\pi_0(\Omega)$ is a totally disconnected compact Hausdorff space and $\pi$ is a closed mapping.
\item If $\Omega$ is metrizable, then $\pi_0(\Omega)$ is metrizable.\qed
\end{enumerate}
\end{theorem}

Theorem~\ref{thm:component-space} allows us to prove the following lemma, which we will use several times.

\begin{lemma}\label{lem:clopen-correspondence}
Let $\Omega$ be a compact Hausdorff space and let $\pi\colon \Omega\to\pi_0(\Omega)$ be the quotient map.
Then $\Lambda\mapsto\pi(\Lambda)$ is a bijection from the family of clopen subsets of $\Omega$ onto the family
of clopen subsets of $\pi_0(\Omega)$, with inverse $\Gamma\mapsto\pi^{-1}(\Gamma)$. This bijection preserves
inclusions, finite unions and intersections, and complements.
\end{lemma}

\begin{proof}
If $\Gamma\subseteq\pi_0(\Omega)$ is clopen, then the preimage $\pi^{-1}(\Gamma)$ is clopen because $\pi$ is continuous. 
Conversely, let $\Lambda\subseteq \Omega$ be clopen. Since $\Lambda$ and $\Omega\setminus\Lambda$ are unions of components 
of $\Omega$, we have $\Lambda=\pi^{-1}\bigl(\pi(\Lambda)\bigr)$ and 
$\pi(\Lambda)\cap\pi(\Omega\setminus\Lambda)=\varnothing$. By Theorem~\ref{thm:component-space}~(i) the mapping $\pi$ is 
closed. Hence both $\pi(\Lambda)$ and $\pi(\Omega\setminus\Lambda)=\pi_0(\Omega)\setminus\pi(\Lambda)$ are closed, which 
means that $\pi(\Lambda)$ is clopen. The two assignments are mutually inverse since $\pi\bigl(\pi^{-1}(\Gamma)\bigr)=\Gamma$. The 
remaining assertions are clear.
\end{proof}

Another lemma that we will need is the following.

\begin{lemma}\label{lem:extend-clopen-partition}
Let $\Omega$ be a totally disconnected compact Hausdorff space. If $\Lambda$ is a closed subset of $\Omega$ such that 
$\Lambda=\Lambda_1\cup\Lambda_2$, where $\Lambda_1$ and $\Lambda_2$ are disjoint clopen subsets of $\Lambda$, then there 
exist disjoint clopen subsets $\Gamma_1$ and $\Gamma_2$ of $\Omega$ such that $\Omega=\Gamma_1\cup\Gamma_2$ and 
$\Gamma_j\cap\Lambda=\Lambda_j$ for $j=1,2$.
\end{lemma}

\begin{proof}
If $\Lambda_1=\varnothing$, then $\Gamma_1=\varnothing$ and $\Gamma_2=\Omega$ have the required properties. Assume
therefore that $\Lambda_1\ne\varnothing$. Recall first that a totally disconnected compact Hausdorff space is zero-dimensional. Since 
$\Lambda$ is closed in $\Omega$ and $\Lambda_1,\Lambda_2$ are closed in $\Lambda$, both are closed, and hence compact, in
$\Omega$. If $x\in\Lambda_1$, then $\Omega\setminus\Lambda_2$ is an open neighbourhood of $x$. Hence, there is a clopen set 
$\Gamma_x$ with $x\in\Gamma_x$ and $\Gamma_x\cap\Lambda_2=\varnothing$. By compactness of $\Lambda_1$, there are points 
$x_1,\ldots,x_k\in\Lambda_1$ such that $\Lambda_1\subseteq\Gamma_{x_1}\cup\cdots\cup\Gamma_{x_k}$. Denote this union by 
$\Gamma_1$. Then $\Gamma_1$ is clopen, $\Lambda_1\subseteq\Gamma_1$ and $\Gamma_1\cap\Lambda_2=\varnothing$, which 
implies $\Gamma_1\cap\Lambda=\Lambda_1$. Now put $\Gamma_2=\Omega\setminus\Gamma_1$, and the statement follows.
\end{proof}

\subsection*{Perfect and scattered spaces}
A topological space $\Omega$ is \emph{perfect} if it has no isolated points, and a subset $\Pi\subseteq\Omega$ is
perfect if it is closed and perfect in its relative topology. A topological space is \emph{scattered} if each of
its non-empty subsets has a point which is isolated in that subset. A space $\Omega$ is called a \emph{Polish
space} if it is separable and completely metrizable. In particular, every compact metrizable space is Polish. We
shall use the Cantor--Bendixson theorem in the form of \cite[Theorem~I.6.4]{Kec95} which says that every Polish
space is the union of a perfect subset and a countable open subset. Together with \cite[Corollary~I.6.5]{Kec95},
this implies that an uncountable Polish space has cardinality $\mathfrak c$ and contains a subset homeomorphic
to the Cantor space $\mathbf 2^{\mathbb N}$ (see the next subsection).

\subsection*{Binary trees}
We consider $\mathbf 2=\{0,1\}$ as a topological space equipped with the discrete topology. For every $k\in\mathbb N$ we denote by 
$\mathbf 2^k$ the set of all $k$-tuples $s=(s_1,\ldots,s_k)$ with $s_j\in\mathbf 2$ and we put $\mathbf 2^0=\{\varnothing\}$, where 
$\varnothing$ means the empty sequence. The set $\mathbf 2^{<\mathbb N}=\bigcup_{k=0}^{\infty}\mathbf 2^k$ is the
\emph{full binary tree}. The elements of $\mathbf 2^{<\mathbb N}$ are \emph{finite binary sequences}. The \emph{length} of 
$s\in\mathbf 2^k$ is $|s|=k$. In particular, $|\varnothing|=0$. If $|s|\ge 1$, then let $s|j=(s_1,\ldots,s_j)$, for $1\le j\le k$, and 
$s|0=\varnothing$. For $t\in \mathbf 2^{<\mathbb N}$, we write $s\subseteq t$ if $s=t|j$ for some $j\le |t|$. In this case $s$ is an 
\emph{initial segment} of $t$. In particular, $\varnothing\subseteq t$ for every $t$. We denote by $s0=(s,0)$ and $s1=(s,1)$ the two 
\emph{immediate successors} of $s$.

The set of all \emph{infinite binary sequences} is denoted by $\mathbf 2^{\mathbb N}$. It is a topological space equipped with the 
product topology. It is a perfect totally disconnected compact metrizable space, called the \emph{Cantor space}. For 
$\omega\in\mathbf 2^{\mathbb N}$, the initial segments of $\omega$ are the finite sequences 
$\omega|j=(\omega_1,\ldots,\omega_j)\in\mathbf 2^j$, where $j\ge 0$. Hence every $\omega\in\mathbf 2^{\mathbb N}$ determines 
the increasing chain $\varnothing=\omega|0\subseteq\omega|1\subseteq\omega|2\subseteq\cdots$ in $\mathbf 2^{<\mathbb N}$.
Let $k\ge 0$ be an integer and $s\in\mathbf 2^{k}$. Then the set $[s]=\{\omega\in\mathbf 2^{\mathbb N};\;\omega|k=s\}$ is clopen in 
$\mathbf 2^{\mathbb N}$. Moreover, these sets form a base of the topology of $\mathbf 2^{\mathbb N}$. We refer the reader to 
\cite[\S I.2]{Kec95} for basics about trees.

Let $\Omega$ be a topological space. A family $\{\Delta_s;\; s\in\mathbf 2^{<\mathbb N}\}$ of non-empty clopen subsets of 
$\Omega$ is a \emph{dyadic family of clopen sets} if $\Delta_\varnothing=\Omega$,
$\Delta_s=\Delta_{s0}\cup\Delta_{s1}$, and $\Delta_{s0}\cap\Delta_{s1}=\varnothing$ for every $s\in\mathbf 2^{<\mathbb N}$. 
The following lemma, whose proof we include for completeness, is a useful tool in what follows.

\begin{lemma}\label{lem:scattered-characterization}
Let $\Omega$ be a totally disconnected compact Hausdorff space. The following assertions are equivalent:
\begin{enumerate}
\item $\Omega$ is not scattered;
\item $\Omega$ contains a non-empty perfect subset;
\item $\Omega$ admits a dyadic family of clopen sets;
\item there exists a continuous surjection of $\Omega$ onto the Cantor space $\mathbf 2^{\mathbb N}$.
\end{enumerate}
If, in addition, $\Omega$ is metrizable, then these assertions are equivalent to
\begin{enumerate}
\setcounter{enumi}{4}
\item $\Omega$ is uncountable.
\end{enumerate}
\end{lemma}

\begin{proof}
We will show that (i) and (ii) are equivalent and that (ii)$\Rightarrow$(iv)$\Rightarrow$(iii)$\Rightarrow$(ii).
Throughout, we use the following immediate consequence of the above definition of a dyadic family obtained by
induction on the length of $t$: in a dyadic family $\{\Delta_s\}$ one has $\Delta_t\subseteq\Delta_s$ whenever
$s\subseteq t$, and for each $k$ the sets $\Delta_t$ with $|t|=k$ are pairwise disjoint with union the whole
space.

(i)$\Leftrightarrow$(ii) Assume (i). Since $\Omega$ is not scattered, there exists a non-empty subset $\Lambda$ of $\Omega$ without
isolated points. Since $\Lambda$ has no isolated points, neither does its closure $\overline{\Lambda}$. Thus $\overline{\Lambda}$ is a
non-empty perfect subset of $\Omega$. The opposite implication is obvious, since a non-empty perfect subset of $\Omega$ has no
isolated points.

(ii)$\Rightarrow$(iv) Let $\Pi$ be a non-empty perfect subset of $\Omega$. Since it is closed in $\Omega$, it is
a totally disconnected compact Hausdorff space and therefore zero-dimensional. Put $\Pi_\varnothing=\Pi$ and
suppose that a non-empty clopen subset $\Pi_s$ of $\Pi$ has been defined. No point of $\Pi_s$ is isolated in
$\Pi_s$, because $\Pi_s$ is open in the perfect space $\Pi$. Hence $\Pi_s$ contains two distinct points $x$ and
$y$. Choose a clopen subset $\Gamma$ of $\Pi$ with $x\in\Gamma$ and $y\notin\Gamma$, and put
$\Pi_{s0}=\Pi_s\cap\Gamma$ and $\Pi_{s1}=\Pi_s\setminus\Gamma$. This yields a dyadic family
$\{\Pi_s;\; s\in\mathbf 2^{<\mathbb N}\}$ of clopen subsets of $\Pi$.

Let $\theta\colon\Pi\to\mathbf 2^{\mathbb N}$ send $x$ to the unique $\omega$ with $x\in\Pi_{\omega|k}$ for all
$k\ge 0$. Recall that the sets $[s]$ form a base of the topology of $\mathbf 2^{\mathbb N}$. Since
$\theta^{-1}([s])=\Pi_s$ (and $\Pi_s$ is open), the mapping $\theta$ is continuous. For every $\omega$ the sets
$\Pi_{\omega|k}$ form a decreasing sequence of non-empty compact sets, whence
$\bigcap_k\Pi_{\omega|k}\neq\varnothing$ and therefore $\theta$ is surjective.

Write $\theta=(\theta_k)_{k\in\mathbb N}$ with $\theta_k\colon\Pi\to\mathbf 2$. In other words,
$\theta_k=\mathrm{pr}_k\circ\,\theta$, where $\mathrm{pr}_k\colon\mathbf 2^{\mathbb N}\to\mathbf 2$ is the
$k$-th coordinate projection. For each $k$ the sets $\theta_k^{-1}(0)$ and $\theta_k^{-1}(1)$ are disjoint clopen subsets of $\Pi$ with 
union $\Pi$. By Lemma~\ref{lem:extend-clopen-partition}, there exists a clopen subset $\Lambda_k$ of $\Omega$ such that
$\Lambda_k\cap\Pi=\theta_k^{-1}(1)$ and $(\Omega\setminus\Lambda_k)\cap\Pi=\theta_k^{-1}(0)$. Define
$\phi\colon\Omega\to\mathbf 2^{\mathbb N}$ by letting the $k$-th coordinate of $\phi(x)$ be $1$ if
$x\in\Lambda_k$ and $0$ otherwise. Since every $\phi^{-1}([s])$ is a finite intersection of clopen sets, the
mapping $\phi$ is continuous, and $\phi|_\Pi=\theta$. Hence
$\phi(\Omega)\supseteq\theta(\Pi)=\mathbf 2^{\mathbb N}$.

(iv)$\Rightarrow$(iii) Let $\phi\colon\Omega\to\mathbf 2^{\mathbb N}$ be a continuous surjection. Then the sets
$\Delta_s=\phi^{-1}([s])$ are clopen and non-empty, $\Delta_\varnothing=\Omega$, and
$[s]=[s0]\cup[s1]$ is a disjoint union, hence so is $\Delta_s=\Delta_{s0}\cup\Delta_{s1}$.

(iii)$\Rightarrow$(ii) Let $\{\Delta_s;\; s\in \mathbf 2^{<\mathbb N}\}$ be a dyadic family of clopen subsets of
$\Omega$ and let $\mathcal F$ be the family of all closed subsets $\Lambda$ of $\Omega$ meeting every
$\Delta_s$. Then $\Omega\in\mathcal F$, and the intersection of a chain in $\mathcal F$ again belongs to
$\mathcal F$, because for each $s$ the corresponding sets $\Lambda\cap\Delta_s$ form a chain of non-empty
compact sets. By Zorn's lemma, $\mathcal F$ has a minimal element $\Pi$, which is non-empty because it meets
$\Delta_\varnothing$.

Note that the subset $\Pi_0$ of all non-isolated points of $\Pi$ is closed. Let $s$ be fixed. For every $k$ the
set $\Delta_s$ contains the $2^k$ pairwise disjoint non-empty sets $\Delta_t$ with $s\subseteq t$
and $|t|=|s|+k$, each of which meets $\Pi$. Hence $\Delta_s\cap\Pi$ is infinite. Being infinite and compact,
it cannot be discrete. This means that it has a point $x$ that is not isolated in it. The point $x$ is not isolated in $\Pi$ either, since a 
point isolated in $\Pi$ would be isolated in $\Delta_s\cap\Pi$. Thus $\Pi_0\in\mathcal F$, and minimality of $\Pi$ gives $\Pi_0=\Pi$.
This shows that $\Pi$ is a non-empty perfect subset of $\Omega$.

Assume finally that $\Omega$ is metrizable. If (iv) holds, then
$\card(\Omega)\ge\card\bigl(\mathbf 2^{\mathbb N}\bigr)=\mathfrak c$ and therefore (v) holds. Conversely,
assume (v) and suppose that $\Omega$ were scattered. Since it is compact and metrizable, $\Omega$ is a Polish
space. By the Cantor--Bendixson theorem \cite[Theorem~I.6.4]{Kec95} we have $\Omega=\Pi\cup\Gamma$ with
$\Pi$ perfect and $\Gamma$ countable. Scatteredness forces $\Pi=\varnothing$, which implies $\Omega=\Gamma$,
that is, $\Omega$ is countable, contradicting (v). Hence (v) implies (i).
\end{proof}
The algebraic counterpart of a dyadic family of clopen sets is a \emph{dyadic family of idempotents}
$\{u_s;\; s\in\mathbf 2^{<\mathbb N}\}$ in $\mathscr A$. It is defined as a family of non-zero idempotents such that
$u_\varnothing=e$, $u_s=u_{s0}+u_{s1}$, and $u_{s0}u_{s1}=0$ for every $s\in\mathbf 2^{<\mathbb N}$.
In the following lemma, we consider the multiplicative structure of such a family.

\begin{lemma} \label{lem:multiplicative-structure-dyadic-family}
For a dyadic family of idempotents $\{u_s;\;s\in\mathbf 2^{<\mathbb N}\}\subseteq \mathscr A$ the following hold:
\begin{enumerate}
\item $u_{s1}u_{s0}=0$ for every $s\in\mathbf 2^{<\mathbb N}$;
\item $u_tu_s=u_su_t=u_t$ whenever $s,t\in\mathbf 2^{<\mathbb N}$ are such that $s\subseteq t$;
\item $u_su_t=0$ whenever $s,t\in\mathbf 2^{<\mathbb N}$ are such that $s\ne t$ and $|s|=|t|$;
\item $u_su_t=u_tu_s$ for arbitrary $s,t\in\mathbf 2^{<\mathbb N}$; moreover, $u_su_t=0$ if
$s\not\subseteq t$ and $t\not\subseteq s$;
\item $\sum_{|s|=m\atop t\subseteq s}u_s=u_t$ for every $m\ge 0$ and for every
$t\in \mathbf 2^{<\mathbb N}$ with $|t|\le m$; in particular, $\sum_{|s|=m}u_s=e$ for every $m\ge 0$;
\item $\sum_{|s|=j\atop s_k=1}u_s=\sum_{|t|=k\atop t_k=1}u_t$ whenever $k,j\in\mathbb N$ are such that
$k\le j$.
\end{enumerate}
\end{lemma}

\begin{proof}
(i) Let $s\in\mathbf 2^{<\mathbb N}$. Since $u_s=u_{s0}+u_{s1}$, we have
$$u_s=u_s^2=(u_{s0}+u_{s1})^2=u_{s0}+u_{s1}+u_{s0}u_{s1}+u_{s1}u_{s0}
=u_s+u_{s0}u_{s1}+u_{s1}u_{s0},$$
which gives $u_{s0}u_{s1}+u_{s1}u_{s0}=0$. Because $u_{s0}u_{s1}=0$ by definition of a dyadic family, we see
that $u_{s1}u_{s0}=0$.

(ii) Let $s\in\mathbf 2^{<\mathbb N}$. Then, by definition, $u_{s0}u_s=u_{s0}(u_{s0}+u_{s1})=u_{s0}$ and
$u_s u_{s0}=(u_{s0}+u_{s1})u_{s0}=u_{s0}$. Similarly, $u_{s1}u_s=u_{s1}$ and $u_s u_{s1}=u_{s1}$.
Suppose that $t\in\mathbf 2^{<\mathbb N}$ is such that $s\subseteq t$. Then a simple induction along the chain
from $s$ to $t$ gives $u_tu_s=u_su_t=u_t$.

(iii) Let $s, t\in\mathbf 2^{<\mathbb N}$ be such that $s\ne t$ and $|s|=|t|$. Denote by $q$ their longest common
initial segment (it can be $\varnothing$). Since $s\ne t$, the length of $q$ is strictly smaller than the length of $s$
and $t$. There is no loss of generality if we assume that $q0\subseteq s$ and $q1\subseteq t$. By (ii),
$u_s=u_s u_{q0}$ and $u_t=u_{q1}u_t$. Hence, $u_s u_t=(u_s u_{q0})(u_{q1} u_t)=u_s(u_{q0}u_{q1})u_t=0$.
Similarly, $u_t u_s=0$.

(iv) Let $s, t\in\mathbf 2^{<\mathbb N}$ be arbitrary. Of course, the assertion holds if $s=t$. Suppose therefore that
$s\ne t$. If $|s|=|t|$, then $u_s u_t=0=u_t u_s$, by (iii). It remains to consider the case when the lengths $|s|$
and $|t|$ are not equal. We may assume that $|s|<|t|$. If $s\subseteq t$, then $u_s u_t=u_t=u_t u_s$, by (ii).
In the case $s\not\subseteq t$, let $q$ be the longest common initial segment of $s$ and $t$.
Again, there is no loss of generality if we assume that $q0\subseteq s$ and $q1\subseteq t$. Now a similar
reasoning as in the proof of (iii) gives $u_s u_t=0=u_t u_s$.

(v) Since $e=u_\varnothing$ and $u_\varnothing=u_0+u_1$, we see that $\sum_{|s|=m}u_s=e$ for $m=0$
and $m=1$. Suppose that $\sum_{|s|=m}u_s=e$ holds for some $m\in\mathbb N$. Since $u_s=u_{s0}+u_{s1}$ and
every $t\in\mathbf 2^{<\mathbb N}$ with $|t|=m+1$ is of the form $t=s0$ or $t=s1$ for some
$s\in\mathbf 2^{<\mathbb N}$ with $|s|=m$, we have
$e=\sum_{|s|=m}u_s=\sum_{|s|=m}(u_{s0}+u_{s1})=\sum_{|t|=m+1}u_t$. This proves
$\sum_{|s|=m}u_s=e$ for every $m\ge 0$. Let $t\in \mathbf 2^{<\mathbb N}$ with $|t|\le m$ be arbitrary.
Multiplying $\sum_{|s|=m}u_s=e$ by $u_t$ gives $\sum_{|s|=m}u_su_t=u_t$. In the sum on the left,
$u_su_t=u_s$ if $t\subseteq s$, by (ii), and $u_su_t=0$ if $t\not\subseteq s$, by (iv). It follows that
$\sum_{|s|=m\atop t\subseteq s}u_s=u_t$.

(vi) Let $k\le j$. If $|s|=j$, then $s|k$ is the only $t\in\mathbf 2^{k}$ with $t\subseteq s$, and the last
entry of $s|k$ is $s_k$. Hence, as $t$ runs through all elements of $\mathbf 2^{k}$ with $t_k=1$, the sets
$\{s\in\mathbf 2^{j};\;t\subseteq s\}$ are pairwise disjoint and their union is
$\{s\in\mathbf 2^{j};\;s_k=1\}$. Applying (v) to each $u_t$ we therefore obtain
$$\sum_{|t|=k\atop t_k=1}u_t=\sum_{|t|=k\atop t_k=1}\ \sum_{|s|=j\atop t\subseteq s}u_s
=\sum_{|s|=j\atop s_k=1}u_s .$$
\end{proof}

We close this subsection with a classical construction of Cantor sets from fast-decreasing sequences of
positive numbers.

\begin{lemma}\label{lem:subsum-cantor-set}
Let $(\tau_k)_{k\in\mathbb N}$ be a sequence of positive numbers such that $\sum_{k>j}\tau_k<\tau_j$ for every
$j\in\mathbb N$. Then the mapping
$$\theta\colon\mathbf 2^{\mathbb N}\to\mathbb R,\qquad
\omega=(\omega_1,\omega_2,\ldots)\mapsto \sum_{k=1}^{\infty}\omega_k\tau_k,$$
is a homeomorphism onto its image $\Sigma=\theta\bigl(\mathbf 2^{\mathbb N}\bigr)$. In particular, $\Sigma$ is a
compact subset of $\mathbb R$ which is homeomorphic to $\mathbf 2^{\mathbb N}$. Moreover, if for
$s\in\mathbf 2^{<\mathbb N}$ we denote by $s000\cdots\in\mathbf 2^{\mathbb N}$ the sequence which extends $s$ by
zeros, then
\begin{equation}\label{eq:branch-approximation}
\bigl|\theta(\omega)-\theta(s000\cdots)\bigr|\le\sum_{k>j}\tau_k
\qquad\bigl(\omega\in\mathbf 2^{\mathbb N},\ j\in\mathbb N,\ s=\omega|j\bigr).
\end{equation}
\end{lemma}

\begin{proof}
The assumption with $j=1$ shows that the series $\sum_{k=1}^{\infty}\tau_k$ of positive numbers converges.
Since $0\le\omega_k\tau_k\le\tau_k$ for every $\omega\in\mathbf 2^{\mathbb N}$ and every $k\in\mathbb N$, the
series defining $\theta$ converges for every $\omega$. Hence $\theta$ is well defined.

Let $\omega\in\mathbf 2^{\mathbb N}$, let $j\in\mathbb N$ and put $s=\omega|j$. Then
$\theta(s000\cdots)=\sum_{k\le j}\omega_k\tau_k$ and therefore
$$\bigl|\theta(\omega)-\theta(s000\cdots)\bigr|=\sum_{k>j}\omega_k\tau_k\le\sum_{k>j}\tau_k,$$
which proves \eqref{eq:branch-approximation}. In particular, $|\theta(\omega)-\theta(\eta)|\le\sum_{k>j}\tau_k$
whenever $\omega|j=\eta|j$. Given $\varepsilon>0$, choose $j\in\mathbb N$ with $\sum_{k>j}\tau_k<\varepsilon$.
Then $\theta$ maps the neighbourhood $[\omega|j]$ of $\omega$ into the $\varepsilon$-neighbourhood of
$\theta(\omega)$. Thus $\theta$ is continuous.

To see that $\theta$ is injective, let $\omega,\eta\in\mathbf 2^{\mathbb N}$ with $\omega\ne\eta$ and let $j$ be
the smallest index with $\omega_j\ne\eta_j$, say $\omega_j=1$ and $\eta_j=0$ (hence, $\omega_k=\eta_k$ for
$k<j$). It follows that
$$\theta(\omega)-\theta(\eta)=\tau_j+\sum_{k>j}(\omega_k-\eta_k)\tau_k\ge\tau_j-\sum_{k>j}\tau_k>0 ,$$
by the assumption on $(\tau_k)_{k\in\mathbb N}$. We conclude that $\theta(\omega)\neq\theta(\eta)$.

Thus $\theta$ is a continuous bijection from the compact space $\mathbf 2^{\mathbb N}$ onto the Hausdorff space
$\Sigma$, hence a homeomorphism. Consequently, $\Sigma$ is a compact subset of $\mathbb R$ which is
homeomorphic to $\mathbf 2^{\mathbb N}$.
\end{proof}

Taking $\tau_k=2\cdot3^{-k}$ in Lemma~\ref{lem:subsum-cantor-set} gives $\sum_{k>j}\tau_k=3^{-j}<\tau_j$
and $\sum_{k}\tau_k=1$. Hence $\Sigma$ is the classical Cantor ternary set. In particular, the Cantor ternary set
is homeomorphic to $\mathbf 2^{\mathbb N}$.

\subsection*{Spectral idempotents and spectral classes}
For $a\in\mathscr A$ we denote by $\sigma(a)$ its spectrum. If $\sigma$ is a clopen subset of $\sigma(a)$, then
the characteristic function $\chi_\sigma$ extends to a holomorphic function on a neighbourhood of $\sigma(a)$ and
the element $e_\sigma=\chi_\sigma(a)$, obtained by the Riesz functional calculus, is the \emph{spectral idempotent}
associated with $\sigma$. Idempotents $u_1,\ldots,u_n$ are \emph{pairwise
orthogonal} if $u_ju_k=0$ for $j\neq k$. If $\sigma_1,\ldots,\sigma_n$ are pairwise disjoint non-empty clopen
subsets of $\sigma(a)$ whose union is $\sigma(a)$, then the associated spectral idempotents are non-zero by
\cite[Proposition~I.7.9]{BD73}, they are pairwise orthogonal and $e=e_{\sigma_1}+\cdots+e_{\sigma_n}$.
Conversely, since $e_{\sigma_1}e_{\sigma_2}=\chi_{\sigma_1}(a)\chi_{\sigma_2}(a)=
(\chi_{\sigma_1}\chi_{\sigma_2})(a)=\chi_{\sigma_1\cap\sigma_2}(a)$,
we have $\sigma_1\cap\sigma_2=\varnothing$ whenever $e_{\sigma_1}e_{\sigma_2}=0$, and
$e_{\sigma_1}+e_{\sigma_2}$ is then the spectral idempotent associated with the
clopen set $\sigma_1\cup\sigma_2$.

In the case of a complex Banach space $X$ and a bounded linear operator $T$ on $X$, the spectral idempotents of
$T$ are the \emph{spectral projections} $P_\sigma=\chi_\sigma(T)$ and their ranges $X_\sigma=P_\sigma X$
are the \emph{spectral subspaces} of $T$. A disjoint clopen decomposition
$\sigma(T)=\sigma_1\cup\cdots\cup\sigma_k$ produces the algebraic and topological direct sum
$X=X_{\sigma_1}\oplus\cdots\oplus X_{\sigma_k}$ into closed invariant subspaces.

We recall the two results of \cite{BK26} that will be used most often. The first one, \cite[Proposition~3.1]{BK26},
states that for $n\in\mathbb N$ the following are equivalent: $\mathcal C_n(\mathscr A)\neq\varnothing$; there exists
$a\in\mathscr A$ such that $\sigma(a)$ has exactly $n$ components; there exist non-zero pairwise orthogonal
idempotents $u_1,\ldots,u_n\in\mathscr A$ with $e=u_1+\cdots+u_n$; the algebra $\mathscr A$ contains a
subalgebra which contains $e$ and is isomorphic to $\mathbb C^n$. The second one, \cite[Theorem~3.4]{BK26},
states that if $(p_j)_{j\in\mathbb N}$ is a sequence of non-zero pairwise orthogonal idempotents in $\mathscr A$,
$w_j=\|p_j\|$ and $(\lambda_j)_{j\in\mathbb N}$ is a sequence of complex numbers with
$\sum_j|\lambda_j|w_j<\infty$, then the series $\sum_j\lambda_jp_j$ converges absolutely and its sum $a$ satisfies
$\sigma(a)=\{0\}\cup\{\lambda_j;\;j\in\mathbb N\}$.

\section{The class $\mathcal C_\infty$ and the spectral depth} \label{sec:Cinfty-spectral-depth}
Recall that
$$\mathcal C_\infty(\mathscr A)=\bigcap_{n=1}^{\infty}\mathcal C_n(\mathscr A)$$
consists precisely of those elements of $\mathscr A$ whose spectra have infinitely many components. The
following proposition is the analogue of \cite[Proposition~3.1]{BK26} for $n=\infty$.

\begin{proposition}\label{prop:Cinfty-idempotent-characterization}
The following assertions are equivalent:
\begin{enumerate}
\item $\mathcal C_\infty(\mathscr A)\neq\varnothing$;
\item the algebra $\mathscr A$ contains a sequence of non-zero pairwise orthogonal idempotents;
\item there exist $a\in\mathscr A$ and a sequence $(\lambda_j)_{j\in\mathbb N}$ of pairwise distinct non-zero
complex numbers such that $\lambda_j\to0$ and $\sigma(a)=\{0\}\cup\{\lambda_j;\;j\in\mathbb N\}$.
\end{enumerate}
\end{proposition}

\begin{proof}
(i)$\Rightarrow$(ii) Choose $a\in\mathcal C_\infty(\mathscr A)$. By Lemma~\ref{lem:infinitely-many-components}
there exists a sequence $(\sigma_j)_{j\in\mathbb N}$ of non-empty pairwise disjoint clopen subsets of $\sigma(a)$.
Let $p_j$ be the spectral idempotent associated with $\sigma_j$. By \cite[Proposition~I.7.9]{BD73}, applied to the
clopen partition $\sigma_1,\ldots,\sigma_j,\sigma(a)\setminus(\sigma_1\cup\cdots\cup\sigma_j)$, each $p_j$ is
non-zero, and $p_jp_k=\chi_{\sigma_j\cap\sigma_k}(a)=0$ for $j\neq k$. Hence $(p_j)_{j\in\mathbb N}$ is a
sequence of non-zero pairwise orthogonal idempotents.

(ii)$\Rightarrow$(iii) This holds by \cite[Corollary~3.5]{BK26}.

(iii)$\Rightarrow$(i) Since the numbers $\lambda_j$ are pairwise distinct, non-zero, and converge to $0$, each of
them is an isolated point of $\sigma(a)$. Hence, the singletons $\{\lambda_j\}$ are pairwise distinct components of $\sigma(a)$.
\end{proof}

\begin{corollary}\label{cor:totally-disconnected-spectrum}
If $a\in\mathscr A$ is such that $\sigma(a)$ is infinite and totally disconnected, then $a\in\mathcal
C_\infty(\mathscr A)$.
\end{corollary}

\begin{proof}
All components of $\sigma(a)$ are singletons, and there are infinitely many of them.
\end{proof}

We now turn to the spectral depth
$$d_\sigma(\mathscr A)=\sup\{n\in\mathbb N;\;\mathcal C_n(\mathscr A)\neq\varnothing\}\in\mathbb N\cup\{\infty\}$$
introduced in \cite{BK26}. If $d_\sigma(\mathscr A)=n$ is finite, then
$$\mathscr A=\mathcal C_1(\mathscr A)\supsetneq\cdots\supsetneq\mathcal C_n(\mathscr A)\supsetneq
\mathcal C_{n+1}(\mathscr A)=\varnothing,$$
by \cite[Corollary~3.3]{BK26}. Proposition~3.1 in \cite{BK26} gives an element whose spectrum has exactly $n$ components, as well as 
the estimate $\dim(\mathscr A)\ge n$. Thus, for finite values, the supremum defining $d_\sigma(\mathscr A)$ is attained. The situation 
for $d_\sigma(\mathscr A)=\infty$ is more delicate. The equality $d_\sigma(\mathscr A)=\infty$ asserts that every class 
$\mathcal C_n(\mathscr A)$ is non-empty, whereas $\mathcal C_\infty(\mathscr A)\neq\varnothing$ asserts that their intersection is 
non-empty. Since a decreasing chain of non-empty sets may have empty intersection, the former does not formally imply the latter.
The question we consider here is whether the supremum defining $d_\sigma(\mathscr A)$ is always attained. 
In Section~\ref{sec:orthogonal-infinitude} it is answered in the negative (Corollary~\ref{cor:infinite-spectral-depth-not-attained}). 
By \cite[Proposition~3.1]{BK26} and Proposition~\ref{prop:Cinfty-idempotent-characterization}, the question is equivalent to a 
statement about idempotents which involves neither spectra nor topology.

\begin{proposition}\label{prop:depth-attained-idempotent-form}
The following assertions are equivalent:
\begin{enumerate}
\item if $d_\sigma(\mathscr A)=\infty$, then $\mathcal C_\infty(\mathscr A)\neq\varnothing$;
\item if $\mathscr A$ contains, for every $n\in\mathbb N$, a family of $n$ non-zero pairwise orthogonal
idempotents, then it contains a sequence of non-zero pairwise orthogonal idempotents;
\item if $\mathscr A$ contains, for every $n\in\mathbb N$, a subalgebra which contains $e$ and is isomorphic to
$\mathbb C^n$, then it contains a sequence of non-zero pairwise orthogonal idempotents.
\end{enumerate}
In particular, the supremum defining $d_\sigma(\mathscr A)$ is attained for every unital complex Banach algebra if and only if (ii) holds 
for every such algebra.
\end{proposition}

\begin{proof}
Fix $n\in\mathbb N$. By \cite[Proposition~3.1]{BK26}, we have $\mathcal C_n(\mathscr A)\neq\varnothing$ if and
only if there exist $n$ non-zero pairwise orthogonal idempotents in $\mathscr A$ summing to $e$, and this holds if
and only if $\mathscr A$ contains a subalgebra which contains $e$ and is isomorphic to $\mathbb C^n$. Moreover,
a family $u_1,\ldots,u_n$ of $n$ non-zero pairwise orthogonal idempotents can always be completed to such a
decomposition of $e$: either $u_1+\cdots+u_n=e$, or $u_1,\ldots,u_n,e-(u_1+\cdots+u_n)$ is a family of $n+1$
non-zero pairwise orthogonal idempotents summing to $e$, whence $\mathcal C_{n+1}(\mathscr A)\neq\varnothing$
and therefore $\mathcal C_n(\mathscr A)\neq\varnothing$. Hence the existence of a family of $n$ non-zero pairwise
orthogonal idempotents, the existence of such a family summing to $e$, and the existence of a subalgebra
containing $e$ and isomorphic to $\mathbb C^n$ all say that $\mathcal C_n(\mathscr A)\neq\varnothing$. Letting
$n$ run through $\mathbb N$, each of the three hypotheses is equivalent to $d_\sigma(\mathscr A)=\infty$.

By Proposition~\ref{prop:Cinfty-idempotent-characterization}, $\mathcal C_\infty(\mathscr A)\neq\varnothing$ if
and only if $\mathscr A$ contains a sequence of non-zero pairwise orthogonal idempotents, which is the
common conclusion of (ii) and (iii).
\end{proof}

The difficulty is that the finite families of idempotents provided by $d_\sigma(\mathscr A)=\infty$ need not be compatible: nothing 
guarantees that a family of $n$ non-zero pairwise orthogonal idempotents can be refined to a family of $n+1$ of them. The following 
observation shows that the implication does hold whenever the spectral classes are not only non-empty but large, which explains
why the phenomenon is invisible in the operator setting.

\begin{proposition}\label{prop:dense-classes}
If $\mathcal C_n(\mathscr A)$ is dense in $\mathscr A$ for every $n\in\mathbb N$, then $\mathcal C_\infty(\mathscr A)$ is a dense 
$G_\delta$-subset of $\mathscr A$. In particular, $\mathcal C_\infty(\mathscr A)\neq\varnothing$ for such algebras.
\end{proposition}

\begin{proof}
By \cite[Theorem~4.2]{BK26}, every class $\mathcal C_n(\mathscr A)$ is an open subset of $\mathscr A$. Being a
Banach space, $\mathscr A$ is a Baire space. Hence, the intersection 
$\mathcal C_\infty(\mathscr A)=\bigcap_{n\in\mathbb N}\mathcal C_n(\mathscr A)$ of countably many dense open 
sets is a dense $G_\delta$-subset of $\mathscr A$.
\end{proof}

\begin{corollary}\label{cor:BX-infinite}
If $X$ is an infinite-dimensional complex Banach space, then $\mathcal C_\infty(\mathscr B(X))$ is non-empty. The same holds 
for every closed subalgebra of $\mathscr B(X)$ which contains $I$ and the ideal of all finite-rank operators in $\mathscr B(X)$.
\end{corollary}

\begin{proof}
For $\mathscr B(X)$ this follows from Proposition~\ref{prop:dense-classes} and \cite[Theorem~4.7]{BK26}. Alternatively, 
\cite[Corollary~3.8]{BK26} produces directly a compact operator lying in $\mathcal C_\infty(\mathscr B(X))$. For the subalgebras, 
apply \cite[Corollary~3.9]{BK26}.
\end{proof}

The converse of Proposition~\ref{prop:dense-classes} is false. Indeed, by \cite[Example~4.9]{BK26} there 
exists a unital commutative complex Banach algebra $\mathscr A$ for which no class $\mathcal C_n(\mathscr A)$ with 
$n\ge2$ is dense in $\mathscr A$, although $\mathcal C_\infty(\mathscr A)\neq\varnothing$. 
In Section~\ref{sec:orthogonal-infinitude} we show that $d_\sigma(\mathscr A)=\infty$ does not imply
$\mathcal C_\infty(\mathscr A)\neq\varnothing$ in general, whereas in Section~\ref{sec:commutative-banach-algebras} we prove 
that it does for every unital commutative Banach algebra.

\section{Orthogonally unbounded and orthogonally infinite algebras}
\label{sec:orthogonal-infinitude}

We will call $\mathscr A$ an \emph{orthogonally unbounded} algebra if for every $n\in\mathbb N$ there exist $n$ non-zero pairwise 
orthogonal idempotents $u_1,\ldots,u_n\in\mathscr A$ such that $e=u_1+\cdots+u_n$. The algebra $\mathscr A$ is 
\emph{orthogonally infinite} if it contains a sequence of non-zero pairwise orthogonal idempotents.
In ring theory, an algebra which is not orthogonally infinite is called \emph{$I$-finite}, while a decomposition of the identity
into $n$ non-zero pairwise orthogonal idempotents is called an \emph{$n$-frame}; see \cite[Section~3]{HoroubNicholson2022}.
As shown in the proof of Proposition~\ref{prop:depth-attained-idempotent-form}, $\mathscr A$ is orthogonally 
unbounded if and only if $d_\sigma(\mathscr A)=\infty$, and it is orthogonally infinite if and only if 
$\mathcal C_\infty(\mathscr A)\neq\varnothing$.

Every orthogonally infinite algebra is orthogonally unbounded. Indeed, if $(u_j)_{j\in\mathbb N}$ is a sequence of non-zero pairwise 
orthogonal idempotents, then for every $n\ge2$ the idempotents $u_1,\ldots,u_{n-1}$ and $e-\sum_{j=1}^{n-1}u_j$ are pairwise 
orthogonal and sum to $e$, and the last one is non-zero because its product with $u_n$ equals $u_n$. Examples of orthogonally 
infinite unital complex Banach algebras are provided by the results of Section~\ref{sec:Cinfty-spectral-depth}. For instance, every 
algebra all of whose spectral classes are dense is orthogonally infinite (Proposition~\ref{prop:dense-classes}). It follows, by 
Corollary~\ref{cor:BX-infinite}, that $\mathscr B(X)$ for an infinite-dimensional complex Banach space $X$, together with each of its 
closed subalgebras containing $I$ and the finite-rank operators, is orthogonally infinite.

The purpose of this section is to show that an orthogonally unbounded algebra need not be orthogonally infinite. We construct
a unital complex Banach algebra which is orthogonally unbounded but $I$-finite. Equivalently, the supremum defining 
$d_\sigma(\mathscr A)$ need not be attained. In Section~\ref{sec:commutative-banach-algebras} we show that this cannot happen 
in the commutative setting: every orthogonally unbounded commutative ring with identity, and in particular every orthogonally 
unbounded unital commutative complex Banach algebra, is orthogonally infinite.

Let $\Omega=\prod_{n=2}^{\infty}\{1,\ldots,n\}$, where each factor carries the discrete topology and $\Omega$ is equipped
with the product topology. By the Tychonoff theorem $\Omega$ is a compact Hausdorff space. We equip $\Omega\times\Omega$ with the product topology. We
first define a map 
$\vartheta\colon \Omega\times \Omega\to \Omega\times \Omega$. If $x=y$, put $\vartheta(x,y)=(x,y)$. Otherwise, let $n$ be the least integer 
such that $x_n\neq y_n$, and obtain $\vartheta(x,y)$ by exchanging the $n$-th coordinates of $x$ and $y$, leaving all other coordinates
unchanged. Let $\vartheta_i\colon \Omega\times \Omega\to \Omega $ $(i=1,2)$ be such that $\vartheta(x,y)=(\vartheta_1(x,y),\vartheta_2(x,y))$.

We claim that $\vartheta$ is continuous. Indeed, fix $k\geq2$ and let $A_k$ consist of all pairs $(x,y)\in \Omega\times \Omega$
such that $x_j=y_j$ for every $2\leq j<k$. Since each factor of $\Omega$ carries the discrete topology,
the set $A_k$ is clopen. By the definition of $\vartheta$, we have
$$ (\vartheta_1(x,y))_k=
\begin{cases}
y_k,&(x,y)\in A_k,\\
x_k,&(x,y)\notin A_k
\end{cases}
\qquad \text{and}\qquad
(\vartheta_2(x,y))_k=
\begin{cases}
x_k,&(x,y)\in A_k,\\
y_k,&(x,y)\notin A_k.
\end{cases}$$
Since the coordinate projections are continuous and $A_k$ is clopen, both functions are continuous.
As this holds for every $k\geq2$, the maps $\vartheta_1$ and $\vartheta_2$, and hence $\vartheta$, are
continuous by the definition of the product topology.

Put $K=[0,1]\times\Omega\times\Omega$, equipped with the product topology. It is a compact Hausdorff space. For $f\in C(\Omega)$, define
$$D_f(x,y)=
\begin{pmatrix}
f(x)&0\\
0&f(y)
\end{pmatrix}$$
for all $x,y\in \Omega$. Let $\mathscr A$ consist of all continuous functions $F\colon K\to \mathbb M_2$ for which there exists 
a function $f\in C(\Omega)$ satisfying 
\begin{equation}\label{eq:matrix-algebra-boundary-conditions}
F(0,x,y)=D_f(x,y) \qquad\text{and}\qquad F(1,x,y)=D_f(\vartheta(x,y))\qquad (x,y\in \Omega).
\end{equation}
The function $f$ associated with $F$ is uniquely determined by $F$, since $f(x)$ is the upper-left entry of $F(0,x,x)$. 
We equip $\mathscr A$ with pointwise operations and the supremum norm $\|F\|=\sup_{(t,x,y)\in K}\|F(t,x,y)\|$,
where the matrix norm is the usual operator norm.

\begin{proposition}\label{prop:matrix-algebra-construction}
The algebra $\mathscr A$ is a unital complex Banach algebra which is orthogonally unbounded.
\end{proposition}

\begin{proof}
If $F,G\in\mathscr A$ correspond to $f,g\in C(\Omega)$ and $\lambda\in\mathbb C$, then $F+G$, $\lambda F$, and $FG$
correspond to $f+g$, $\lambda f$, and $fg$, respectively. Moreover, the constant function $e=I_2$ belongs to
$\mathscr A$. If $F_j\in\mathscr A$ converge uniformly to $F$, then the associated functions $f_j$ converge
uniformly to the function $f$, where $f(x)$ is the upper-left entry of $F(0,x,x)$. Passing to the limit in
\eqref{eq:matrix-algebra-boundary-conditions} gives $F\in\mathscr A$. Hence $\mathscr A$ is a closed unital
subalgebra of $C(K,\mathbb M_2)$.

Fix $n\geq2$. For $1\le r\le n$, let $f_{n,r}$ be the characteristic function of the clopen set
$\{x\in\Omega;\;x_n=r\}$. Let $V_n$ be the clopen set of all pairs $(x,y)$ whose first differing coordinate is
the $n$-th one. Put
$$R(t)=\begin{pmatrix}1-t&-t\\ t&1-t\end{pmatrix}\qquad(t\in[0,1]).$$
Since $\det R(t)=(1-t)^2+t^2>0$, the function $U_n$, defined by $U_n(t,x,y)=R(t)$ if $(x,y)\in V_n$ and by
$U_n(t,x,y)=I_2$ otherwise, is continuous with invertible values. Let
$P_{n,r}=U_n D_{f_{n,r}}U_n^{-1}$. Then $P_{n,1},\ldots,P_{n,n}$ are non-zero pairwise orthogonal idempotents
with sum $e$, because the same holds for $D_{f_{n,1}},\ldots,D_{f_{n,n}}$.

It remains to show that $P_{n,r}\in\mathscr A$. Since $R(0)=I_2$, the first condition in
\eqref{eq:matrix-algebra-boundary-conditions} holds. Note that $f_{n,r}(x)$ depends only on $x_n$. If
$(x,y)\in V_n$, then $\vartheta$ exchanges exactly the $n$-th coordinates of $x$ and $y$, and therefore
$D_{f_{n,r}}(\vartheta(x,y))$ is obtained from $D_{f_{n,r}}(x,y)$ by exchanging the diagonal entries. Conjugation
by $R(1)=\left(\begin{smallmatrix}0&-1\\ 1&0\end{smallmatrix}\right)$ does exactly this. If $(x,y)\notin V_n$, then
$\vartheta$ does not change the $n$-th coordinates, and therefore $D_{f_{n,r}}(\vartheta(x,y))=D_{f_{n,r}}(x,y)$. In both cases 
$P_{n,r}(1,x,y)=D_{f_{n,r}}(\vartheta(x,y))$. Thus $P_{n,r}\in\mathscr A$, and $\mathscr A$ is orthogonally unbounded.
\end{proof}

We next examine the functions associated with non-zero idempotents of $\mathscr A$.

\begin{lemma}\label{lem:matrix-algebra-idempotent-boundary}
Let $F\in\mathscr A$ be a non-zero idempotent, and let $f\in C(\Omega)$ be its associated function.
Then $f$ is non-zero, takes values in $\{0,1\}$, and satisfies
\begin{equation}\label{eq:matrix-algebra-swap-identity}
f(x)+f(y)=f(\vartheta_1(x,y))+f(\vartheta_2(x,y))\qquad(x,y\in \Omega).
\end{equation}
\end{lemma}

\begin{proof}
Since $F^2=F$, the first condition in \eqref{eq:matrix-algebra-boundary-conditions} gives $D_f^2=D_f$. Hence $f(x)\in\{0,1\}$ 
for every $x\in\Omega$. For fixed $x,y\in\Omega$, the matrices $F(t,x,y)$, $t\in[0,1]$, are idempotents. The trace of an
idempotent matrix equals its rank. Hence $t\mapsto\tr F(t,x,y)$ is a continuous integer-valued function on $[0,1]$,
and therefore it is constant. Comparing its values at $t=0$ and $t=1$ gives \eqref{eq:matrix-algebra-swap-identity}.
If $f=0$, then $F(0,x,y)=0$, and constancy of the rank gives $F(t,x,y)=0$ for all $(t,x,y)\in K$. This contradicts
$F\ne0$.
\end{proof}

\begin{lemma}\label{lem:matrix-algebra-coordinate-rigidity}
If $f\in C(\Omega)$ is a non-constant $\{0,1\}$-valued function satisfying \eqref{eq:matrix-algebra-swap-identity},
then there exists $n\ge2$ such that $f(x)$ depends only on $x_n$.
\end{lemma}

\begin{proof}
Choose $x,y\in\Omega$ with $f(x)\ne f(y)$. Since $f$ is continuous at $y$, there is $N\ge2$ such that $f(y')=f(y)$ for
the point $y'$ with $y'_k=y_k$ for $k\le N$ and $y'_k=x_k$ for $k>N$. The points $x$ and $y'$ differ in finitely many
coordinates. Changing these coordinates one at a time, we find two points which differ in exactly one coordinate and
have different values of $f$. Let $n$ be the least index for which such a pair exists. Then $f$ does not change if a
single coordinate with index less than $n$ is changed.

We write each $x\in\Omega$ as $x=(p,q,a)$, where $p=(x_2,\ldots,x_{n-1})$, $q=x_n$, and $a=(x_{n+1},x_{n+2},\ldots)$.
We call $a$ the tail of $x$. If $n=2$, then $p$ is empty. Any two choices of $p$ can be connected by changing their
finitely many coordinates one at a time. Hence $f(p,q,a)$ does not depend on $p$, and we may write $f(p,q,a)=g(q,a)$.
If $q\ne q'$, then the first coordinate in which $(p,q,a)$ and $(p,q',b)$ differ is the $n$-th one. By the definition
of $\vartheta$, $\vartheta\bigl((p,q,a),(p,q',b)\bigr)=\bigl((p,q',a),(p,q,b)\bigr)$.
Therefore \eqref{eq:matrix-algebra-swap-identity} gives
\begin{equation}\label{eq:matrix-algebra-g-identity}
g(q,a)+g(q',b)=g(q',a)+g(q,b)
\end{equation}
for all $q,q'\in\{1,\ldots,n\}$ and all tails $a,b$. If $q=q'$, then \eqref{eq:matrix-algebra-g-identity} holds trivially.

By the choice of $n$, there are two points which differ only in the $n$-th coordinate and have different values of $f$. They have the 
same tail $c$. Hence there are $r,s\in\{1,\ldots,n\}$ such that $g(r,c)=0$ and $g(s,c)=1$. Let $b$ be an arbitrary tail. 
By \eqref{eq:matrix-algebra-g-identity}, $g(r,b)+g(s,c)=g(s,b)+g(r,c)$, that is, $g(r,b)+1=g(s,b)$. Since $g$ takes values in $\{0,1\}$, we 
obtain $g(r,b)=0$. Now let $q\in\{1,\ldots,n\}$ and let $a,b$ be arbitrary tails. By \eqref{eq:matrix-algebra-g-identity} with $q'=r$, 
we have $g(q,a)-g(q,b)=g(r,a)-g(r,b)=0$. Thus $g(q,a)$ does not depend on $a$, and therefore $f(x)$ depends only on $x_n$.
\end{proof}

\begin{proposition}\label{prop:matrix-algebra-I-finite}
The algebra $\mathscr A$ is not orthogonally infinite.
\end{proposition}

\begin{proof}
Suppose that $(F_j)_{j\in\mathbb N}$ is a sequence of non-zero pairwise orthogonal idempotents in $\mathscr A$, and
let $f_j$ be the function associated with $F_j$. By the first condition in \eqref{eq:matrix-algebra-boundary-conditions}, 
$f_jf_k=0$ for $j\ne k$. By Lemma~\ref{lem:matrix-algebra-idempotent-boundary}, every $f_j$ is non-zero. Hence no $f_j$ is constant. 
By Lemma~\ref{lem:matrix-algebra-coordinate-rigidity}, there exist $n_j\ge2$ and a non-empty set
$S_j\subseteq\{1,\ldots,n_j\}$ such that $f_j(x)=1$ if and only if $x_{n_j}\in S_j$. If $n_j\ne n_k$, then there exists
$x\in\Omega$ with $x_{n_j}\in S_j$ and $x_{n_k}\in S_k$, which contradicts $f_jf_k=0$. Hence all $n_j$ are equal to
some $n$, and the sets $S_j$ are pairwise disjoint non-empty subsets of $\{1,\ldots,n\}$. This is impossible.
Thus $\mathscr A$ contains no sequence of non-zero pairwise orthogonal idempotents.
\end{proof}

\begin{corollary}\label{cor:infinite-spectral-depth-not-attained}
There exists a unital complex Banach algebra $\mathscr A$ such that $\mathcal C_n(\mathscr A)\neq\varnothing$ for
every $n\in\mathbb N$, whereas $\mathcal C_\infty(\mathscr A)=\varnothing$. In particular, $d_\sigma(\mathscr A)=\infty$.
\end{corollary}

\begin{proof}
Let $\mathscr A$ be the algebra constructed above. By Proposition~\ref{prop:matrix-algebra-construction}, $\mathscr A$
is orthogonally unbounded. Hence $\mathcal C_n(\mathscr A)\neq\varnothing$ for every $n\in\mathbb N$, by
\cite[Proposition~3.1]{BK26}. By Proposition~\ref{prop:matrix-algebra-I-finite}, $\mathscr A$ is not
orthogonally infinite. Hence $\mathcal C_\infty(\mathscr A)=\varnothing$, by
Proposition~\ref{prop:Cinfty-idempotent-characterization}.
\end{proof}

\section{The cardinality of the set of spectral components} \label{sec:cardinality-spectral-components}
Since the spectrum of an element of a Banach algebra is a non-empty compact subset of $\mathbb C$, the set of its components has 
cardinality at most $\mathfrak c$. At first sight it seems plausible that, depending on the Continuum Hypothesis, intermediate 
cardinalities might occur. The following theorem shows that within ZFC this is not the case.

\begin{theorem}\label{thm:possibility-for-cardinality}
For every $a\in\mathscr A$ we have
$$\card\bigl(\pi_0(\sigma(a))\bigr)\in\mathbb N\cup\{\aleph_0,\mathfrak c\}.$$
\end{theorem}

\begin{proof}
The spectrum $\sigma(a)$ is a non-empty compact metrizable space and therefore, by Theorem~\ref{thm:component-space}, 
the space $\pi_0(\sigma(a))$ is a totally disconnected compact metrizable space. In particular, it is Polish. By the
Cantor--Bendixson theorem \cite[Theorem~I.6.4]{Kec95} there exist a perfect subset $\Pi$ and a countable open subset $\Gamma$ of 
$\pi_0(\sigma(a))$ with $\pi_0(\sigma(a))=\Pi\cup\Gamma$. If $\Pi=\varnothing$, then $\pi_0(\sigma(a))$ is countable, and therefore 
its cardinality lies in $\mathbb N\cup\{\aleph_0\}$. Otherwise, $\Pi$, being a non-empty perfect subset of a Polish space, has cardinality 
$\mathfrak c$ by \cite[Corollary~I.6.5]{Kec95}. Since $\card\bigl(\pi_0(\sigma(a))\bigr)\le\card\bigl(\sigma(a)\bigr)\le\mathfrak c$, 
we conclude that $\card\bigl(\pi_0(\sigma(a))\bigr)=\mathfrak c$.
\end{proof}

Accordingly, for $\nu\in\{\aleph_0,\mathfrak c\}$ we put
$$\mathcal C_{\nu}(\mathscr A)=\{a\in\mathscr A;\;\card\bigl(\pi_0(\sigma(a))\bigr)=\nu\},$$
keeping the earlier meaning of $\mathcal C_n(\mathscr A)$ for $n\in\mathbb N$. The classes 
$\mathcal C_{\aleph_0}(\mathscr A)$ and $\mathcal C_{\mathfrak c}(\mathscr A)$ are disjoint and, by
Theorem~\ref{thm:possibility-for-cardinality}, 
$\mathcal C_\infty(\mathscr A)=\mathcal C_{\aleph_0}(\mathscr A)\cup\mathcal C_{\mathfrak c}(\mathscr A)$.
Moreover, by Proposition~\ref{prop:Cinfty-idempotent-characterization}, if 
$\mathcal C_\infty(\mathscr A)\neq\varnothing$, then already $\mathcal C_{\aleph_0}(\mathscr A)\neq\varnothing$.

Theorem~\ref{thm:possibility-for-cardinality} also allows us to refine the spectral depth into a cardinal-valued invariant. From now on, 
we put
\begin{equation}\label{eq:depth-definition}
d_\sigma(\mathscr A)=\sup\{\card\bigl(\pi_0(\sigma(a))\bigr);\;a\in\mathscr A\},
\end{equation}
the supremum being taken in the ordered set of cardinals. By Theorem~\ref{thm:possibility-for-cardinality}, we have 
$d_\sigma(\mathscr A)\in\mathbb N\cup\{\aleph_0,\mathfrak c\}$. This is consistent with the definition used in \cite{BK26}. Indeed, by 
\cite[Proposition~3.1]{BK26} we have $\mathcal C_n(\mathscr A)\neq\varnothing$ if and only if some element has a spectrum with 
exactly $n$ components. Hence \eqref{eq:depth-definition} agrees with $\sup\{n;\;\mathcal C_n(\mathscr A)\neq\varnothing\}$ 
whenever the latter is finite, and $d_\sigma(\mathscr A)\in\{\aleph_0,\mathfrak c\}$ if and only if 
$\mathcal C_n(\mathscr A)\neq\varnothing$ for every $n\in\mathbb N$. In particular, if $d_\sigma(\mathscr A)$ is finite, then the 
supremum in \eqref{eq:depth-definition} is attained. If $d_\sigma(\mathscr A)=\mathfrak c$, then 
$\mathcal C_{\mathfrak c}(\mathscr A)\neq\varnothing$, because a supremum of cardinals from $\mathbb N\cup\{\aleph_0\}$ never 
exceeds $\aleph_0$. Hence the supremum in \eqref{eq:depth-definition} can fail to be attained only when 
$d_\sigma(\mathscr A)=\aleph_0$. This indeed happens: by Corollary~\ref{cor:infinite-spectral-depth-not-attained}, the algebra 
constructed in Section~\ref{sec:orthogonal-infinitude} satisfies $d_\sigma(\mathscr A)=\aleph_0$ and 
$\mathcal C_{\aleph_0}(\mathscr A)=\varnothing$.\medskip

The classes $\mathcal C_\nu(\mathscr A)$ are stable under the holomorphic functional calculus in the following sense.

\begin{proposition}\label{prop:functional-calculus}
Let $a\in\mathscr A$ and let $f$ be a holomorphic function on an open neighbourhood $U$ of $\sigma(a)$. Then
\begin{enumerate}
\item $\card\bigl(\pi_0(\sigma(f(a)))\bigr)\le\card\bigl(\pi_0(\sigma(a))\bigr)$;
\item if $f$ is injective on $U$, then $\card\bigl(\pi_0(\sigma(f(a)))\bigr)=\card\bigl(\pi_0(\sigma(a))\bigr)$;
consequently, for every $\nu\in\mathbb N\cup\{\aleph_0,\mathfrak c\}$ we have $a\in\mathcal C_\nu(\mathscr A)$
if and only if $f(a)\in\mathcal C_\nu(\mathscr A)$.
\end{enumerate}
\end{proposition}

\begin{proof}
(i) By the spectral mapping theorem \cite[Theorem~VII.4.10]{Con90} we have $\sigma(f(a))=f(\sigma(a))$. The
assertion follows from Lemma~\ref{lem:continuous-image-components} applied to the continuous mapping
$f|_{\sigma(a)}$.

(ii) If $f$ is injective and holomorphic on $U$, then by the open mapping theorem \cite[Theorem~7.5]{Con78} and
\cite[Corollary~7.6]{Con78} the mapping $f\colon U\to f(U)$ is biholomorphic. Hence the restriction
$f|_{\sigma(a)}\colon\sigma(a)\to f(\sigma(a))=\sigma(f(a))$ is a homeomorphism. Homeomorphic spaces have the
same number of components.
\end{proof}

It follows from Proposition~\ref{prop:functional-calculus} that the classes $\mathcal C_\nu(\mathscr A)$, 
$\nu\in\mathbb N\cup\{\aleph_0,\mathfrak c\}$, are invariant under affine transformations and under taking the inverse. Indeed, if 
$a\in\mathscr A$ and $\alpha,\beta\in\mathbb C$, $\alpha\ne 0$, then $a\in\mathcal C_\nu(\mathscr A)$ if and only if 
$\alpha a+\beta e\in\mathcal C_\nu(\mathscr A)$. If $a$ is invertible, then $a\in\mathcal C_\nu(\mathscr A)$ if and only if 
$a^{-1}\in\mathcal C_\nu(\mathscr A)$. Since an algebra isomorphism preserves spectra, the spectral depth is an isomorphism 
invariant.\medskip

We now characterize the largest of the three regimes. The following theorem is the exact analogue, one level
higher, of \cite[Theorem~3.4]{BK26}. As we shall see, a dyadic family of idempotents, which is a purely algebraic 
concept, produces a single element whose spectrum is as disconnected as a compact subset of $\mathbb C$ can be.

\begin{theorem}\label{thm:dyadic-idempotents-Cantor-spectrum}
Suppose that $\{u_s;\;s\in\mathbf 2^{<\mathbb N}\}$ is a dyadic family of idempotents in $\mathscr A$. 
Then there exists $a\in\mathscr A$ such that $\sigma(a)$ is homeomorphic to the Cantor space
$\mathbf 2^{\mathbb N}$. In particular, 
$\mathcal C_{\mathfrak c}(\mathscr A)\neq\varnothing$ and $d_\sigma(\mathscr A)=\mathfrak c$.
\end{theorem}

\begin{proof}
For every $k\in\mathbb N$ let $M_k=\sum_{|s|=k}\|u_s\|$ and let $p_k=\sum_{|s|=k,\,s_k=1}u_s$. The element $p_k$ is an idempotent, 
by Lemma~\ref{lem:multiplicative-structure-dyadic-family}~(iii), and $\|p_k\|\le M_k$. Put $\tau_0=1$ and 
$\tau_k=\min\bigl\{\frac{1}{3}\tau_{k-1},\ \frac{1}{2^{k}M_k}\bigr\}$ for $k\in\mathbb N$. Then $\tau_k\le \frac{1}{3^{k-j}}\tau_j$ 
for $k>j$ and therefore $\sum_{k>j}\tau_k\le\tau_j\sum_{k>j}\frac{1}{3^{k-j}}\le \frac{1}{2}\tau_j<\tau_j$.
Lemma~\ref{lem:subsum-cantor-set} provides a homeomorphism $\theta$ of $\mathbf 2^{\mathbb N}$ onto the compact
set $\Sigma=\theta\bigl(\mathbf 2^{\mathbb N}\bigr)$. Since $\tau_k\|p_k\|\le2^{-k}$, the series
$a=\sum_{k=1}^{\infty}\tau_kp_k$ converges absolutely. Its partial sums $c_j=\sum_{k=1}^{j}\tau_kp_k$ satisfy
$$\|a-c_j\|\le\sum_{k>j}\tau_k\|p_k\|\le2^{-j}\qquad(j\in\mathbb N).$$

Fix $j\in\mathbb N$. By Lemma~\ref{lem:multiplicative-structure-dyadic-family}~(vi), we have
$p_k=\sum_{|s|=j,\,s_k=1}u_s$ for $k\le j$. Hence, by rearranging the finite sums,
$$c_j=\sum_{k=1}^{j}\tau_k\sum_{\substack{|s|=j\\ s_k=1}}u_s=\sum_{|s|=j}\Bigl(\sum_{k=1}^{j}s_k\tau_k\Bigr)u_s
=\sum_{|s|=j}\lambda_su_s,$$
where $\lambda_s=\sum_{k=1}^{j}s_k\tau_k=\theta(s000\cdots)$.
We claim that $\sigma(c_j)=\{\lambda_s;\;|s|=j\}$. Recall that the idempotents $u_s$ with $|s|=j$ are non-zero by definition, and that 
they are pairwise orthogonal and sum to $e$ by Lemma~\ref{lem:multiplicative-structure-dyadic-family}~(iii) and~(v).
If $\zeta\ne\lambda_s$ for every $s$ with $|s|=j$, then $c_j-\zeta e=\sum_{|s|=j}(\lambda_s-\zeta)u_s$. 
Since $\bigl(\sum_{|s|=j}(\lambda_s-\zeta)u_s\bigr)\bigl(\sum_{|s|=j}(\lambda_s-\zeta)^{-1}u_s\bigr)=\sum_{|s|=j}u_s=e$, by
the orthogonality, and the two factors commute, we see that $c_j-\zeta e$ is invertible, that is, $\zeta\notin\sigma(c_j)$.
On the other hand, by the orthogonality, $(c_j-\lambda_se)u_s=0$ and therefore $c_j-\lambda_s e$ is not invertible since $u_s\ne 0$. 
This proves the claim; in particular, $\sigma(c_j)\subseteq\Sigma$.

Since the idempotents $u_s$ commute with each other, by Lemma~\ref{lem:multiplicative-structure-dyadic-family}~(iv),
the element $c_j$ commutes with $a-c_j$. Recall that $\sigma(x+y)\subseteq\sigma(x)+\sigma(y)$ for commuting
$x,y\in\mathscr A$ \cite[Theorem~11.23]{Rud91}. Note that the spectra $\sigma(\pm(a-c_j))$ are contained in the closed disc $D_j$
with centre $0$ and radius $\|a-c_j\|\le2^{-j}$. Hence
$$\sigma(a)\subseteq\sigma(c_j)+D_j\subseteq\Sigma+D_j\qquad\text{and}\qquad\sigma(c_j)\subseteq\sigma(a)+D_j
\qquad(j\in\mathbb N).$$
The first inclusion gives $\dist(\zeta,\Sigma)\le2^{-j}$ for every $\zeta\in\sigma(a)$ and every $j$. Since $\Sigma$ is
closed, we conclude that $\sigma(a)\subseteq\Sigma$. Conversely, let $\zeta=\theta(\omega)\in\Sigma$. The second inclusion gives
$\dist\bigl(\lambda_{\omega|j},\sigma(a)\bigr)\le2^{-j}$, while $\lambda_{\omega|j}\to\zeta$ by
\eqref{eq:branch-approximation}. Since $\sigma(a)$ is closed, we have $\zeta\in\sigma(a)$. Hence $\sigma(a)=\Sigma$.

By Lemma~\ref{lem:subsum-cantor-set}, the spectrum $\sigma(a)=\Sigma$ is homeomorphic to
$\mathbf 2^{\mathbb N}$. In particular, it is totally disconnected of cardinality $\mathfrak c$. Hence
$\card\bigl(\pi_0(\sigma(a))\bigr)=\mathfrak c$, that is, $a\in\mathcal C_{\mathfrak c}(\mathscr A)$ and
$d_\sigma(\mathscr A)=\mathfrak c$.
\end{proof}

\begin{corollary}\label{cor:continuum-spectral-class-tree-characterization}
The following assertions are equivalent:
\begin{enumerate}
\item $\mathcal C_{\mathfrak c}(\mathscr A)\neq\varnothing$;
\item there exist $a\in\mathscr A$ and a dyadic family $\{\sigma_s;\;s\in\mathbf 2^{<\mathbb N}\}$ of clopen
subsets of $\sigma(a)$;
\item $\mathscr A$ contains a dyadic family of idempotents;
\item there exists $a\in\mathscr A$ such that $\sigma(a)$ is homeomorphic to the Cantor space $\mathbf 2^{\mathbb N}$.
\end{enumerate}
\end{corollary}

\begin{proof}
(i)$\Rightarrow$(ii) Choose $a\in\mathcal C_{\mathfrak c}(\mathscr A)$ and let
$\pi\colon\sigma(a)\to\pi_0(\sigma(a))$ be the quotient map. By Theorem~\ref{thm:component-space} the space
$\Omega=\pi_0(\sigma(a))$ is a totally disconnected compact metrizable space, and by assumption it is
uncountable. Hence Lemma~\ref{lem:scattered-characterization} provides a dyadic family $\{\Delta_s;\;s\in\mathbf
2^{<\mathbb N}\}$ of clopen subsets of $\Omega$. Put $\sigma_s=\pi^{-1}(\Delta_s)$. By
Lemma~\ref{lem:clopen-correspondence} each $\sigma_s$ is a non-empty clopen subset of $\sigma(a)$, and the
relations $\sigma_\varnothing=\sigma(a)$, $\sigma_s=\sigma_{s0}\cup\sigma_{s1}$ and
$\sigma_{s0}\cap\sigma_{s1}=\varnothing$ are preserved.

(ii)$\Rightarrow$(iii) For $s\in\mathbf 2^{<\mathbb N}$ let $u_s$ be the spectral idempotent of $a$ associated
with $\sigma_s$. By \cite[Proposition~I.7.9]{BD73} each $u_s$ is non-zero, $u_\varnothing=\chi_{\sigma(a)}(a)=e$,
and, as recalled in Section~\ref{sec:preliminaries}, $u_{s0}u_{s1}=\chi_{\sigma_{s0}\cap\sigma_{s1}}(a)=0$ and
$u_{s0}+u_{s1}=\chi_{\sigma_{s0}\cup\sigma_{s1}}(a)=\chi_{\sigma_s}(a)=u_s$.

(iii)$\Rightarrow$(iv) This is Theorem~\ref{thm:dyadic-idempotents-Cantor-spectrum}.

(iv)$\Rightarrow$(i) The Cantor space $\mathbf 2^{\mathbb N}$ is totally disconnected of cardinality $\mathfrak c$.
\end{proof}

Combining the results of this section with Proposition~\ref{prop:Cinfty-idempotent-characterization} we obtain a
description of the three regimes of the spectral depth purely in terms of the idempotents of $\mathscr A$.

\begin{corollary}\label{cor:depth-trichotomy}
Exactly one of the following holds.
\begin{enumerate}
\item There exists $n\in\mathbb N$ such that every family of non-zero pairwise orthogonal idempotents in $\mathscr A$
has at most $n$ members and some such family has $n$ members. In this case $d_\sigma(\mathscr A)=n$.
\item $\mathscr A$ contains families of non-zero pairwise orthogonal idempotents of arbitrarily large finite
cardinality, but no dyadic family of idempotents. In this case $d_\sigma(\mathscr A)=\aleph_0$, and $\mathcal
C_{\aleph_0}(\mathscr A)\neq\varnothing$ if and only if $\mathscr A$ contains a sequence of non-zero
pairwise orthogonal idempotents.
\item $\mathscr A$ contains a dyadic family of idempotents. In this case $d_\sigma(\mathscr A)=\mathfrak c$ and
$\mathcal C_{\mathfrak c}(\mathscr A)\neq\varnothing$.
\end{enumerate}
\end{corollary}

\begin{proof}
If $\mathscr A$ contains a dyadic family of idempotents, then $d_\sigma(\mathscr A)=\mathfrak c$ by
Theorem~\ref{thm:dyadic-idempotents-Cantor-spectrum}, which gives (iii). Suppose that $\mathscr A$ contains no
dyadic family of idempotents. Then $\mathcal C_{\mathfrak c}(\mathscr A)=\varnothing$ by
Corollary~\ref{cor:continuum-spectral-class-tree-characterization} and therefore $d_\sigma(\mathscr A)\le\aleph_0$. If the
cardinalities of families of non-zero pairwise orthogonal idempotents are bounded, let $n$ be the largest
cardinality occurring. Then $\mathcal C_n(\mathscr A)\neq\varnothing$ and $\mathcal C_{n+1}(\mathscr
A)=\varnothing$ by \cite[Proposition~3.1]{BK26} and the completion argument in the proof of
Proposition~\ref{prop:depth-attained-idempotent-form}. Hence, $d_\sigma(\mathscr A)=n$. Otherwise 
$\mathcal C_n(\mathscr A)\neq\varnothing$ for every $n$, which gives $d_\sigma(\mathscr A)=\aleph_0$. The last assertion is
Proposition~\ref{prop:Cinfty-idempotent-characterization}.
\end{proof}

Both alternatives in case (ii) of Corollary~\ref{cor:depth-trichotomy} really occur. By
Corollary~\ref{cor:infinite-spectral-depth-not-attained}, an algebra of spectral depth $\aleph_0$ need not contain a sequence of
non-zero pairwise orthogonal idempotents. On the other hand, $\mathscr B(X_{\mathrm{AH}})$ (Example~\ref{ex:AH}) and the algebra
$c$ of all convergent complex sequences (Example~\ref{ex:commutative-examples}~(ii)) have spectral depth $\aleph_0$ and contain 
such a sequence. More generally, by Proposition~\ref{prop:dense-classes} and Theorem~\ref{thm:commutative-depth}, an algebra in 
case (ii) contains such a sequence whenever its spectral classes are dense or it is commutative. 

We now record the operator-theoretic reformulation of Corollary~\ref{cor:continuum-spectral-class-tree-characterization}. Recall that 
the spectral idempotents of an operator are its spectral projections and that their ranges are the corresponding spectral subspaces.

\begin{corollary}\label{cor:continuum-spectral-class-operator-characterization}
Let $X$ be a complex Banach space. The following assertions are equivalent:
\begin{enumerate}
\item $\mathcal C_{\mathfrak c}(\mathscr B(X))\neq\varnothing$;
\item there exist $T\in\mathscr B(X)$ and a family $\{X_s;\;s\in\mathbf 2^{<\mathbb N}\}$ of non-zero closed
spectral subspaces of $T$ such that $X_\varnothing=X$ and $X_s=X_{s0}\oplus X_{s1}$ for $s\in\mathbf 2^{<\mathbb N}$;
\item there is a family $\{P_s;\;s\in\mathbf 2^{<\mathbb N}\}$ of non-zero bounded projections on $X$ such that
$P_\varnothing=I$ and $P_s=P_{s0}+P_{s1}$, $P_{s0}P_{s1}=P_{s1}P_{s0}=0$ for $s\in\mathbf 2^{<\mathbb N}$.
\end{enumerate}
\end{corollary}

\begin{proof}
(i)$\Rightarrow$(ii) By the implication (i)$\Rightarrow$(ii) of
Corollary~\ref{cor:continuum-spectral-class-tree-characterization}, there are $T\in\mathscr B(X)$ and a dyadic family 
$\{\sigma_s;\;s\in\mathbf 2^{<\mathbb N}\}$ of clopen subsets of $\sigma(T)$. Put $X_s=P_{\sigma_s}X$. These spectral subspaces 
are closed, and they are non-zero by \cite[Proposition~I.7.9]{BD73}, since $\sigma_s\neq\varnothing$. Moreover, $X_\varnothing=X$, 
and $\sigma_s=\sigma_{s0}\cup\sigma_{s1}$ with $\sigma_{s0}\cap\sigma_{s1}=\varnothing$ gives $X_s=X_{s0}\oplus X_{s1}$.

(ii)$\Rightarrow$(iii) For each $s$ choose a clopen subset $\sigma_s$ of $\sigma(T)$ with $X_s=P_{\sigma_s}X$ and put 
$P_s=P_{\sigma_s}$. Then $P_s\neq0$ because $X_s\neq\{0\}$, and $P_\varnothing=I$ because $X_\varnothing=X$. Since the spectral 
projections of $T$ commute, the range of $P_{s0}P_{s1}=P_{s1}P_{s0}$ is contained in $X_{s0}\cap X_{s1}=\{0\}$ and therefore 
$P_{s0}P_{s1}=P_{s1}P_{s0}=0$. Hence $P_{s0}+P_{s1}$ is an idempotent with range $X_{s0}\oplus X_{s1}=X_s$, and it commutes with 
$P_s$. Two commuting idempotents with the same range coincide; hence $P_{s0}+P_{s1}=P_s$.

(iii)$\Rightarrow$(i) Since $\{P_s;\;s\in\mathbf 2^{<\mathbb N}\}$ is a dyadic family of idempotents in $\mathscr B(X)$, (i) follows from 
the implication (iii)$\Rightarrow$(i) of Corollary~\ref{cor:continuum-spectral-class-tree-characterization}.
\end{proof}

In other words, as soon as some operator on $X$ has spectrum with $\mathfrak c$ components, the space $X$ admits a full binary tree 
of topological direct decompositions, all of whose members are spectral subspaces of a single operator. We next exhibit a large class of 
spaces where this happens.

\begin{proposition}\label{prop:unconditional-basis}
Let $X$ be a complex Banach space which contains a complemented closed subspace $Y$ with an unconditional
Schauder basis. Then $d_\sigma(\mathscr B(X))=\mathfrak c$ and $\mathcal C_\nu(\mathscr B(X))\neq\varnothing$ for
every $\nu\in\mathbb N\cup\{\aleph_0,\mathfrak c\}$.
\end{proposition}

\begin{proof}
Let $R$ be a bounded projection of $X$ onto $Y$. Let $(e_j)_{j\in\mathbb N}$ be an unconditional Schauder basis of
$Y$ and let $(f_j)_{j\in\mathbb N}$ be the corresponding biorthogonal functionals. For $\Lambda\subseteq\mathbb N$
we define $Q_\Lambda x=\sum_{j\in\Lambda}f_j(Rx)\,e_j$ $(x\in X)$. Since the basis is unconditional, $Q_\Lambda$
is a bounded projection on $X$ by \cite[Proposition~1.c.7]{LT77}. It is obvious that $Q_\Lambda\neq0$ if
$\Lambda\ne\varnothing$ and that $Q_\Lambda(I-R)=(I-R)Q_\Lambda=0$. Moreover, if $\Lambda,\Lambda'\subseteq\mathbb N$
are disjoint, then $Q_\Lambda Q_{\Lambda'}=0$ and $Q_\Lambda+Q_{\Lambda'}=Q_{\Lambda\cup\Lambda'}$.

For $s\in\mathbf 2^{k}$ let $r_s=\sum_{i=1}^{k}s_i2^{i-1}$ and $\Lambda_s=\{j\in\mathbb N;\;j-1\equiv r_s
\pmod{2^k}\}$. Then $\Lambda_\varnothing=\mathbb N$ and each $\Lambda_s$ is non-empty. Since $r_{s0}=r_s$ and
$r_{s1}=r_s+2^k$, the set $\Lambda_s$ is the disjoint union of $\Lambda_{s0}$ and $\Lambda_{s1}$. If $s$ consists
of zeros only (including $s=\varnothing$), let $P_s=Q_{\Lambda_s}+(I-R)$. Otherwise, let $P_s=Q_{\Lambda_s}$. Note
that $s0$ consists of zeros only if and only if $s$ does, and that $s1$ never does. It follows that
$\{P_s;\;s\in\mathbf 2^{<\mathbb N}\}$ is a family of non-zero bounded projections such that
$P_\varnothing=Q_{\mathbb N}+(I-R)=I$, $P_s=P_{s0}+P_{s1}$ and $P_{s0}P_{s1}=P_{s1}P_{s0}=0$. By
Corollary~\ref{cor:continuum-spectral-class-operator-characterization}, we have
$\mathcal C_{\mathfrak c}(\mathscr B(X))\neq\varnothing$. Hence $d_\sigma(\mathscr B(X))=\mathfrak c$.

Finally, $\mathcal C_{\mathfrak c}(\mathscr B(X))\subseteq\mathcal C_\infty(\mathscr B(X))\subseteq\mathcal
C_n(\mathscr B(X))$ for every $n\in\mathbb N$. In particular, $\mathcal C_\infty(\mathscr B(X))$ is non-empty.
By Proposition~\ref{prop:Cinfty-idempotent-characterization}, there is an operator whose spectrum is of the form
$\{0\}\cup\{\lambda_j;\;j\in\mathbb N\}$, where $\lambda_j\to0$. This operator belongs to
$\mathcal C_{\aleph_0}(\mathscr B(X))$.
\end{proof}

\begin{example}\label{ex:AH}
The Argyros--Haydon space $X_{\mathrm{AH}}$ of \cite{AH11} is an infinite-dimensional separable complex Banach space with the 
property that every operator on $X_{\mathrm{AH}}$ is of the form $\alpha I+K$, where $\alpha\in\mathbb C$ and $K$ is compact. 
Since the spectrum of a compact operator is countable, the spectrum of every operator on $X_{\mathrm{AH}}$ is countable. 
Hence $\mathcal C_{\mathfrak c}\bigl(\mathscr B(X_{\mathrm{AH}})\bigr)=\varnothing$. On the other hand,
$\mathcal C_\infty\bigl(\mathscr B(X_{\mathrm{AH}})\bigr)\neq\varnothing$ by \cite[Corollary~3.8]{BK26}. It follows that
$d_\sigma\bigl(\mathscr B(X_{\mathrm{AH}})\bigr)=\aleph_0$.
By \cite[Theorem~1.1]{KL17}, the space $X_{\mathrm{AH}}$ has a shrinking Schauder basis. Thus the hypothesis of
Proposition~\ref{prop:unconditional-basis} cannot be weakened to the mere existence of a Schauder basis. Moreover,
the equivalent conditions of Corollary~\ref{cor:continuum-spectral-class-operator-characterization} fail for
$X=X_{\mathrm{AH}}$.\qed
\end{example}

Together with Proposition~\ref{prop:unconditional-basis}, Example~\ref{ex:AH} shows that both infinite values of
the spectral depth are attained by algebras of the form $\mathscr B(X)$. In particular, $d_\sigma(\mathscr B(X))$
is a non-trivial isomorphic invariant of the space $X$.

We close the section with the behaviour of the spectral depth under finite direct sums.

\begin{theorem}\label{thm:depth-additive}
Let $\mathscr A$ and $\mathscr B$ be unital complex Banach algebras. Then
$$d_\sigma(\mathscr A\oplus\mathscr B)=d_\sigma(\mathscr A)+d_\sigma(\mathscr B).$$
\end{theorem}

\begin{proof}
Denote by $e$ and $f$ the units of $\mathscr A$ and $\mathscr B$, respectively. Put $\mu=d_\sigma(\mathscr A)$ and 
$\nu=d_\sigma(\mathscr B)$. Recall that $\sigma(c\oplus d)=\sigma(c)\cup\sigma(d)$ for $c\in\mathscr A$ and $d\in\mathscr B$.

For every $c\in\mathscr A$ and $d\in\mathscr B$, inequality \eqref{eq:cardinality-components-inequality} of
Lemma~\ref{lem:union-components} and the monotonicity of cardinal addition give
$$\card\bigl(\pi_0(\sigma(c\oplus d))\bigr)\le\card\bigl(\pi_0(\sigma(c))\bigr)+\card\bigl(\pi_0(\sigma(d))\bigr)
\le\mu+\nu. $$
Hence $d_\sigma(\mathscr A\oplus\mathscr B)\le\mu+\nu$. For the reverse inequality, let $a\in\mathscr A$ and $b\in\mathscr B$ be 
arbitrary. Since $\sigma(a)$ and $\sigma(b)$ are compact, there exists $\lambda\in\mathbb C$ such that 
$\sigma(a)\cap\sigma(b+\lambda f)=\varnothing$. The
set $\sigma(b+\lambda f)=\lambda+\sigma(b)$ has the same number of components as $\sigma(b)$. By
\eqref{eq:cardinality-components-equality},
$$d_\sigma(\mathscr A\oplus\mathscr B)\ge\card\bigl(\pi_0(\sigma(a\oplus(b+\lambda f)))\bigr)
=\card\bigl(\pi_0(\sigma(a))\bigr)+\card\bigl(\pi_0(\sigma(b))\bigr).$$
Taking the supremum over $a\in\mathscr A$ and $b\in\mathscr B$, it remains to observe that
$$\sup_{a,b}\bigl(\card\bigl(\pi_0(\sigma(a))\bigr)+\card\bigl(\pi_0(\sigma(b))\bigr)\bigr)=\mu+\nu.$$
If both $\mu$ and $\nu$ are finite, the suprema are attained by \cite[Proposition~3.1]{BK26} and the equality is
clear. If $\max\{\mu,\nu\}=\mathfrak c$, say $\mu=\mathfrak c$, then, as observed after
\eqref{eq:depth-definition}, there is $a\in\mathcal C_{\mathfrak c}(\mathscr A)$ and the left-hand side is at
least $\mathfrak c=\mu+\nu$. Finally, if $\mu,\nu\le\aleph_0$ and at least one of them equals $\aleph_0$, say
$\mu=\aleph_0$, then $\card\bigl(\pi_0(\sigma(a))\bigr)$ takes arbitrarily large finite values. Hence the
left-hand side is at least $\aleph_0=\mu+\nu$.
\end{proof}

Theorem~\ref{thm:depth-additive} yields the following trichotomy for a finite direct sum
$\mathscr A=\mathscr A_1\oplus\cdots\oplus\mathscr A_m$ of unital complex Banach algebras. First,
$d_\sigma(\mathscr A)=\mathfrak c$ if and only if $d_\sigma(\mathscr A_j)=\mathfrak c$ for at least one $j$. Equivalently,
$\mathcal C_{\mathfrak c}(\mathscr A)\neq\varnothing$ if and only if $\mathcal C_{\mathfrak c}(\mathscr A_j)\neq\varnothing$
for some $j$. Second, $\mathscr A$ has finite spectral depth if and only if all summands do, and then the depths add.
Third, $d_\sigma(\mathscr A)=\aleph_0$ if and only if no summand has spectral depth $\mathfrak c$ and at least
one has infinite spectral depth.

Moreover, $\mathcal C_\infty(\mathscr A)\neq\varnothing$ if and only if $\mathcal C_\infty(\mathscr A_j)\neq\varnothing$ for some 
$j$, and the same holds for $\mathcal C_{\aleph_0}$. This follows from Proposition~\ref{prop:Cinfty-idempotent-characterization}, since
a sequence of non-zero pairwise orthogonal idempotents in $\mathscr A$ has, for some $j$, infinitely many non-zero coordinates
in $\mathscr A_j$.

\section{Commutative Banach algebras} \label{sec:commutative-banach-algebras}
In this section $\mathscr A$ is a commutative complex Banach algebra with unit $e$. We denote by $\Delta$ its character 
space, that is, the set of all non-zero multiplicative linear functionals on $\mathscr A$, equipped with the relative weak-$^*$ topology. 
It is a compact Hausdorff space, and the Gelfand transform $a\mapsto\hat a$, where $\hat a(\phi)=\phi(a)$ for $\phi\in\Delta$, is a 
norm-decreasing unital homomorphism of $\mathscr A$ into $C(\Delta)$ satisfying
\begin{equation}\label{eq:gelfand-spectrum}
\sigma(a)=\hat a(\Delta)\qquad (a\in\mathscr A).
\end{equation}
The kernel of the Gelfand transform is the radical $\rad\mathscr A$, which consists of the quasinilpotent elements of $\mathscr A$,
that is, of elements $x\in\mathscr A$ with spectral radius $r(x)=0$. We refer the reader to \cite{Kan09} for these facts. Our starting point 
is the \v{S}ilov idempotent theorem, for which we refer to \cite[\S 3.5]{Kan09} and \cite[Theorem~21.5]{BD73}. The theorem says that 
for a clopen set $\Gamma\subseteq\Delta$ there exists an idempotent $u\in\mathscr A$ such that $\hat u=\chi_\Gamma$.

\begin{proposition}\label{prop:idempotents-clopen-bijection}
The Gelfand transform induces a bijection $u\mapsto\Gamma_u$ from the set of all idempotents of $\mathscr A$ onto the family of 
all clopen subsets of $\Delta$, determined by $\hat u=\chi_{\Gamma_u}$. Moreover, for idempotents 
$u,v\in\mathscr A$:
\begin{enumerate}
\item $u=0$ if and only if $\Gamma_u=\varnothing$, and $u=e$ if and only if $\Gamma_u=\Delta$;
\item $uv=0$ if and only if $\Gamma_u\cap\Gamma_v=\varnothing$, and in that case $u+v$ is an idempotent with
$\Gamma_{u+v}=\Gamma_u\cup\Gamma_v$.
\end{enumerate}
\end{proposition}

\begin{proof}
If $u$ is an idempotent, then $\hat u$ is an idempotent of $C(\Delta)$, hence $\hat u=\chi_{\Gamma}$ for a uniquely determined 
clopen subset $\Gamma=\Gamma_u$ of $\Delta$. Surjectivity of $u\mapsto\Gamma_u$ follows by the \v{S}ilov idempotent theorem.
For injectivity, let $u,v$ be idempotents with $\hat u=\hat v$ and put $x=u-v$. Since $\mathscr A$ is commutative, we have 
$x^3=u^3-3u^2v+3uv^2-v^3=u-v=x$. Moreover, $\hat x=0$. Hence $x\in\rad\mathscr A$. In particular, $e-x^2$ is invertible. Since 
$x(e-x^2)=x-x^3=0$, it follows that $x=0$, and therefore $u=v$.

(i) If $\Gamma_u=\varnothing$, then $\hat u=0$. By the injectivity just proved, $u=0$. The statement 
about $e$ follows in the same way.

(ii) If $uv=0$, then $\chi_{\Gamma_u\cap\Gamma_v}=\hat u\hat v=0$ and therefore $\Gamma_u\cap\Gamma_v=\varnothing$.
Conversely, if $\Gamma_u\cap\Gamma_v=\varnothing$, then $uv$ is an idempotent with $\widehat{uv}=0$. Hence
$uv=0$ by (i). In that case $(u+v)^2=u+v$ and
$\widehat{u+v}=\chi_{\Gamma_u}+\chi_{\Gamma_v}=\chi_{\Gamma_u\cup\Gamma_v}$.
\end{proof}

Combining Proposition~\ref{prop:idempotents-clopen-bijection} with Lemma~\ref{lem:clopen-correspondence} we obtain a bijection 
between the idempotents of $\mathscr A$ and the clopen subsets of the space $\pi_0(\Delta)$, which by 
Theorem~\ref{thm:component-space} is a totally disconnected compact Hausdorff space. This bijection turns families of non-zero 
pairwise orthogonal idempotents into families of non-empty pairwise disjoint clopen subsets of $\pi_0(\Delta)$, and dyadic families of 
idempotents into dyadic families of clopen subsets of $\pi_0(\Delta)$. This is the mechanism behind the following theorem, which 
computes the spectral depth of an arbitrary unital commutative complex Banach algebra.

\begin{theorem}\label{thm:commutative-depth}
The following assertions hold:
\begin{enumerate}
\item for $n\in\mathbb N$ we have $\mathcal C_n(\mathscr A)\neq\varnothing$ if and only if $\card\bigl(\pi_0(\Delta)\bigr)\ge n$;
\item $\mathcal C_\infty(\mathscr A)\neq\varnothing$ if and only if $\pi_0(\Delta)$ is infinite;
\item $\mathcal C_{\mathfrak c}(\mathscr A)\neq\varnothing$ if and only if $\pi_0(\Delta)$ is not scattered.
\end{enumerate}
Consequently,
$$ d_\sigma(\mathscr A)=
\begin{cases}
\card\bigl(\pi_0(\Delta)\bigr), & \text{if $\pi_0(\Delta)$ is finite},\\[2pt]
\aleph_0, & \text{if $\pi_0(\Delta)$ is infinite and scattered},\\[2pt]
\mathfrak c, & \text{if $\pi_0(\Delta)$ is not scattered},
\end{cases} $$
and in each case the supremum in \eqref{eq:depth-definition} is attained.
\end{theorem}

\begin{proof}
(i) Assume that $\card\bigl(\pi_0(\Delta)\bigr)\ge n$, that is, that $\Delta$ has at least $n$ components. By
\cite[Lemma~2.1]{BK26} there are non-empty pairwise disjoint clopen subsets $\Gamma_1,\ldots,\Gamma_n$ of $\Delta$ with union 
$\Delta$. Let $u_1,\ldots,u_n$ be the corresponding idempotents furnished by Proposition~\ref{prop:idempotents-clopen-bijection}. 
They are non-zero and pairwise orthogonal, and $u_1+\cdots+u_n$ is an idempotent whose Gelfand transform is $\chi_\Delta$. 
Hence $u_1+\cdots+u_n=e$. By \cite[Proposition~3.1]{BK26} we get $\mathcal C_n(\mathscr A)\neq\varnothing$.

Conversely, let $a\in\mathcal C_n(\mathscr A)$. By \eqref{eq:gelfand-spectrum} the mapping 
$\hat a\colon\Delta\to\sigma(a)$ is a continuous surjection. Hence, Lemma~\ref{lem:continuous-image-components} gives
$n\le\card\bigl(\pi_0(\sigma(a))\bigr)\le\card\bigl(\pi_0(\Delta)\bigr)$.

(ii) If $\pi_0(\Delta)$ is infinite, then by Lemma~\ref{lem:infinitely-many-components} there is a sequence
$(\Gamma_j)_{j\in\mathbb N}$ of non-empty pairwise disjoint clopen subsets of $\Delta$, and
Proposition~\ref{prop:idempotents-clopen-bijection} yields a sequence of non-zero pairwise orthogonal idempotents
in $\mathscr A$. Proposition~\ref{prop:Cinfty-idempotent-characterization} gives $\mathcal C_\infty(\mathscr
A)\neq\varnothing$. The converse follows from (i).

(iii) Suppose that $\pi_0(\Delta)$ is not scattered. Let $\pi\colon\Delta\to\pi_0(\Delta)$ be the quotient map. By 
Lemma~\ref{lem:scattered-characterization} there is a dyadic family $\{\Lambda_s;\;s\in\mathbf 2^{<\mathbb N}\}$ of clopen subsets 
of $\pi_0(\Delta)$, and by Lemma~\ref{lem:clopen-correspondence} the sets $\Gamma_s=\pi^{-1}(\Lambda_s)$ form a dyadic family 
of clopen subsets of $\Delta$. By Proposition~\ref{prop:idempotents-clopen-bijection} the corresponding idempotents $u_s$ are 
non-zero and satisfy $u_\varnothing=e$, $u_{s0}u_{s1}=0$ and $u_{s0}+u_{s1}=u_s$, that is, they form a dyadic family of idempotents 
in $\mathscr A$. Theorem~\ref{thm:dyadic-idempotents-Cantor-spectrum} gives 
$\mathcal C_{\mathfrak c}(\mathscr A)\neq\varnothing$. In fact, some element of $\mathscr A$ has spectrum homeomorphic to
the Cantor space $\mathbf 2^{\mathbb N}$.

Conversely, let $\mathcal C_{\mathfrak c}(\mathscr A)\neq\varnothing$. By
Corollary~\ref{cor:continuum-spectral-class-tree-characterization} the algebra $\mathscr A$ contains a dyadic
family $\{u_s;\;s\in\mathbf 2^{<\mathbb N}\}$ of idempotents. By
Proposition~\ref{prop:idempotents-clopen-bijection} the sets $\Gamma_{u_s}$ form a dyadic family of clopen
subsets of $\Delta$, and by Lemma~\ref{lem:clopen-correspondence} the sets $\pi(\Gamma_{u_s})$ form a dyadic
family of clopen subsets of $\pi_0(\Delta)$. By Lemma~\ref{lem:scattered-characterization} the space $\pi_0(\Delta)$ is not
scattered.

Now we will compute the spectral depth. If $\pi_0(\Delta)$ is finite with $n$ elements, then by (i) we have 
$\mathcal C_n(\mathscr A)\neq\varnothing$ and $\mathcal C_{n+1}(\mathscr A)=\varnothing$. Hence $d_\sigma(\mathscr A)=n$ and
the supremum is attained by \cite[Proposition~3.1]{BK26}. If $\pi_0(\Delta)$ is infinite and scattered, then 
$\mathcal C_\infty(\mathscr A)\neq\varnothing$ by (ii) and $\mathcal C_{\mathfrak c}(\mathscr A)=\varnothing$ by (iii). It follows that
$\mathcal C_{\aleph_0}(\mathscr A)\neq\varnothing$ and $d_\sigma(\mathscr A)=\aleph_0$ is attained. If $\pi_0(\Delta)$ is not 
scattered, then $\mathcal C_{\mathfrak c}(\mathscr A)\neq\varnothing$ by (iii), which means that $d_\sigma(\mathscr A)=\mathfrak c$ 
is attained.
\end{proof}

It follows immediately from Theorem~\ref{thm:commutative-depth} that the supremum defining the spectral depth is 
attained for every unital commutative complex Banach algebra. Equivalently, by
Proposition~\ref{prop:depth-attained-idempotent-form}, a unital commutative complex Banach algebra which is orthogonally 
unbounded is orthogonally infinite. By Corollary~\ref{cor:infinite-spectral-depth-not-attained}, commutativity cannot be
omitted here.

Moreover, the result is of a purely algebraic nature: neither the norm nor the completeness of $\mathscr A$ plays any role, and the 
\v{S}ilov idempotent theorem is not needed. To see this, let $\mathscr R$ be a commutative ring with unit $e$ and, for an idempotent
$f\in\mathscr R$, denote by $\rho(f)\in\mathbb N\cup\{\infty\}$ the supremum of those $m\in\mathbb N$ for which
there exist non-zero pairwise orthogonal idempotents $r_1,\ldots,r_m\in\mathscr R$ with
$f=r_1+\cdots+r_m$; we put $\rho(0)=0$.

\begin{proposition}\label{prop:commutative-algebraic}
If $\mathscr R$ contains, for every $n\in\mathbb N$, a family of $n$ non-zero pairwise orthogonal idempotents, then it contains a 
sequence of non-zero pairwise orthogonal idempotents.
\end{proposition}

\begin{proof}
We first observe that $\rho(e)=\infty$. Indeed, let $n\in\mathbb N$ and let $r_1,\ldots,r_n\in\mathscr R$ be
non-zero pairwise orthogonal idempotents. Then $r=r_1+\cdots+r_n$ is an idempotent, and so is $e-r$. If
$e-r\ne0$, then $e=r_1+\cdots+r_n+(e-r)$ is a decomposition of $e$ into $n+1$ non-zero pairwise orthogonal
idempotents; otherwise $e=r_1+\cdots+r_n$. In either case $\rho(e)\ge n$.

Next we claim that if $f\in\mathscr R$ is an idempotent with $\rho(f)=\infty$, then there exist non-zero
orthogonal idempotents $p,q\in\mathscr R$ such that $f=p+q$ and $\rho(q)=\infty$. Since $\rho(f)\ge2$, there
certainly are non-zero orthogonal idempotents $p,q$ with $f=p+q$, and it suffices to show that $\rho(p)$ and
$\rho(q)$ cannot both be finite. Assume that they are, and let $f=r_1+\cdots+r_m$ be an arbitrary decomposition
of $f$ into non-zero pairwise orthogonal idempotents. As $\mathscr R$ is commutative, each $pr_j$ is an
idempotent, $(pr_j)(pr_k)=p^2r_jr_k=0$ for $j\ne k$, and $\sum_{j=1}^{m}pr_j=pf=p$. Hence the non-zero elements
among $pr_1,\ldots,pr_m$ form a decomposition of $p$ into pairwise orthogonal non-zero idempotents, and there are
at most $\rho(p)$ of them. Similarly, at most $\rho(q)$ of the elements $qr_1,\ldots,qr_m$ are non-zero. Since
$r_j=fr_j=pr_j+qr_j$ and $r_j\ne0$, at least one of $pr_j$, $qr_j$ is non-zero for every $j$, and therefore
$m\le\rho(p)+\rho(q)$. As the decomposition of $f$ was arbitrary, $\rho(f)\le\rho(p)+\rho(q)<\infty$, which is
a contradiction.

Put $f_0=e$. Applying the claim recursively, we obtain non-zero idempotents $p_j$ and $f_j$ $(j\in\mathbb N)$
such that $f_{j-1}=p_j+f_j$, $p_jf_j=0$ and $\rho(f_j)=\infty$ for every $j\in\mathbb N$. From
$f_{j}f_{j-1}=f_j$ we get $f_jf_k=f_j$ whenever $k\le j$, and $p_jf_{j-1}=p_j$ gives $p_jf_k=p_j$ for every
$k<j$. Hence, if $k<j$, then $p_kp_j=p_k(f_kp_j)=(p_kf_k)p_j=0$. Thus $(p_j)_{j\in\mathbb N}$ is a sequence of non-zero pairwise 
orthogonal idempotents.
\end{proof}

Together with Proposition~\ref{prop:depth-attained-idempotent-form}, this shows again (independently of 
Theorem~\ref{thm:commutative-depth}) that every unital commutative complex Banach algebra which is orthogonally unbounded 
satisfies $\mathcal C_\infty(\mathscr A)\neq\varnothing$. Commutativity enters only through the fact that $pr_j$ is an idempotent. This 
step is unavailable in the general case, and Corollary~\ref{cor:infinite-spectral-depth-not-attained} shows that commutativity 
cannot be dispensed with.\medskip

Since the character space of $C(\Omega)$ is homeomorphic to $\Omega$ for every compact Hausdorff space $\Omega$,
Theorem~\ref{thm:commutative-depth} contains and completes \cite[Proposition~4.5]{BK26}.

\begin{corollary}\label{cor:continuous-functions-depth}
Let $\Omega$ be a compact Hausdorff space. Then
$$d_\sigma\bigl(C(\Omega)\bigr)=
\begin{cases}
\card\bigl(\pi_0(\Omega)\bigr), & \text{if $\pi_0(\Omega)$ is finite},\\[2pt]
\aleph_0, & \text{if $\pi_0(\Omega)$ is infinite and scattered},\\[2pt]
\mathfrak c, & \text{if $\pi_0(\Omega)$ is not scattered}.
\end{cases}$$
For every unital commutative complex Banach algebra $\mathscr A$ with character space $\Delta$,
$$d_\sigma(\mathscr A)=d_\sigma\bigl(C(\Delta)\bigr)=d_\sigma(\mathscr A/\rad\mathscr A).$$
\end{corollary}

\begin{proof}
The first assertion is Theorem~\ref{thm:commutative-depth} applied to $C(\Omega)$. For the second one, note that the algebras
$\mathscr A$, $C(\Delta)$, and $\mathscr A/\rad\mathscr A$ have homeomorphic character spaces, and apply 
Theorem~\ref{thm:commutative-depth} to each of the three algebras.
\end{proof}

In particular, the spectral depth of a commutative Banach algebra does not see its radical, and it depends only
on the topology of the space of components of its character space.

\begin{example}\label{ex:commutative-examples}
(i) Let $\mathscr A$ be the disc algebra $A(\mathbb D)$ ($\mathbb D$ is the open unit disc in the complex plane), or the algebra 
$C[0,1]$, or any unital commutative complex Banach algebra with connected character space. Then $d_\sigma(\mathscr A)=1$, that is, 
every element has connected spectrum.

(ii) For $\mathscr A=c$, the algebra of all convergent complex sequences, the character space is $\mathbb N_\infty$, which is infinite, 
totally disconnected and scattered. Hence $d_\sigma(c)=\aleph_0$.

(iii) Let $\Omega$ be the Cantor ternary set. Then $\pi_0(\Omega)=\Omega$ is perfect, hence not scattered, and therefore 
$d_\sigma\bigl(C(\Omega)\bigr)=\mathfrak c$.

(iv) Let $\mathbb N$ carry the discrete topology, put $\Omega=\mathbb N\times\overline{\mathbb D}$, and let 
$\mathscr A=C_b(\Omega)$ be the unital commutative complex Banach algebra of all bounded continuous complex functions on 
$\Omega$. Every element of $\mathscr A$ is a bounded sequence $f=(f_j)_{j\in\mathbb N}$ of continuous functions on 
$\overline{\mathbb D}$. An easy computation gives $\sigma(f)=\overline{\bigcup_{j\in\mathbb N}f_j(\overline{\mathbb D})}$. Let 
$\Lambda$ be the Cantor ternary set, let $(\lambda_j)_{j\in\mathbb N}$ be a dense sequence in $\Lambda$ and let $f_j$ be the 
constant function $\lambda_j$. Then $\sigma(f)=\overline{\{\lambda_j;\;j\in\mathbb N\}}=\Lambda$, which shows that 
$f\in\mathcal C_{\mathfrak c}(\mathscr A)$ and $d_\sigma\bigl(C_b(\Omega)\bigr)=\mathfrak c$. This is worth comparing with
\cite[Example~4.9]{BK26}, where it is shown that none of the classes $\mathcal C_n(\mathscr A)$ with $n\ge2$ is dense in 
$\mathscr A$: the spectral classes of an algebra of maximal spectral depth may all fail to be dense.\qed
\end{example}

We conclude the section with the finite case, where the spectral depth is reflected in a decomposition of the algebra itself. Recall that 
$\mathscr A$ is said to have no non-trivial idempotents if its only idempotents are $0$ and its unit $e$.

\begin{theorem}\label{thm:finite-depth-decomposition}
Let $n\in\mathbb N$. The following assertions are equivalent:
\begin{enumerate}
\item $d_\sigma(\mathscr A)=n$;
\item there exist non-zero closed ideals $\mathscr A_1,\ldots,\mathscr A_n$ of $\mathscr A$ such that
$\mathscr A=\mathscr A_1\oplus\cdots\oplus\mathscr A_n$, each $\mathscr A_j$ is a Banach algebra with a unit, and 
$\mathcal C_2(\mathscr A_j)=\varnothing$ for $j=1,\ldots,n$;
\item there exist non-zero pairwise orthogonal idempotents $u_1,\ldots,u_n\in\mathscr A$ with $e=u_1+\cdots+u_n$
such that no unital complex Banach algebra $u_j\mathscr A$ has a non-trivial idempotent.
\end{enumerate}
If these conditions hold, then every decomposition as in (ii) is of the form $\mathscr A_j=u_j\mathscr A$ with $u_1,\ldots,u_n$ as in (iii), 
and it is unique up to a permutation of the summands.
\end{theorem}

\begin{proof}
(i)$\Rightarrow$(iii) and (i)$\Rightarrow$(ii). Let $d_\sigma(\mathscr A)=n$. By \cite[Proposition~3.1]{BK26} there exist non-zero 
pairwise orthogonal idempotents $u_1,\ldots,u_n\in\mathscr A$ such that $e=u_1+\cdots+u_n$. Put $\mathscr A_j=u_j\mathscr A$. 
Since $\mathscr A$ is commutative and $\mathscr A_j$ is the range of the bounded projection $a\mapsto u_ja$, each $\mathscr A_j$ is 
a closed ideal of $\mathscr A$ and a Banach algebra with unit $u_j$. Since $a=ea=\sum_{j=1}^nu_ja$ for every $a\in\mathscr A$, we 
have $\mathscr A=\mathscr A_1+\cdots+\mathscr A_n$. If $a_1+\cdots+a_n=0$ with $a_j\in\mathscr A_j$, then multiplication by $u_k$ 
gives $a_k=u_ka_k=u_k\sum_ja_j=0$, which proves that the sum is direct.

Suppose that $\mathcal C_2(\mathscr A_j)\neq\varnothing$ for some $j$. By \cite[Proposition~3.1]{BK26} there exist non-zero 
orthogonal idempotents $v,w\in\mathscr A_j$ such that $u_j=v+w$. Replacing $u_j$ by $v$ and $w$ in the family $u_1,\ldots,u_n$ we 
obtain $n+1$ non-zero pairwise orthogonal idempotents of $\mathscr A$ summing to $e$, whence $d_\sigma(\mathscr A)\ge n+1$, 
a contradiction. Thus $\mathcal C_2(\mathscr A_j)=\varnothing$ for every $j$, and by \cite[Corollary~3.2]{BK26} the algebra 
$\mathscr A_j$ has no non-trivial idempotent.

(ii)$\Rightarrow$(i) Let $u_j$ be the unit of $\mathscr A_j$. Since the sum is direct and 
$\mathscr A_j\mathscr A_k\subseteq\mathscr A_j\cap\mathscr A_k=\{0\}$ for $j\neq k$, the idempotents $u_1,\ldots,u_n$ are non-zero,
pairwise orthogonal and sum to $e$. By \cite[Proposition~3.1]{BK26} we get 
$\mathcal C_n(\mathscr A)\neq\varnothing$, which gives $d_\sigma(\mathscr A)\ge n$. Conversely, let $v\in\mathscr A$ be an 
idempotent. Then $v=\sum_ju_jv$ and each $u_jv$ is an idempotent of $\mathscr A_j$. Since 
$\mathcal C_2(\mathscr A_j)=\varnothing$, the algebra $\mathscr A_j$ has no non-trivial idempotent by \cite[Corollary~3.2]{BK26}.
It follows that $u_jv\in\{0,u_j\}$. Hence every idempotent of $\mathscr A$ is of the form $\sum_{j\in F}u_j$ for some subset
$F$ of $\{1,\ldots,n\}$. It follows that every family of non-zero pairwise orthogonal idempotents in $\mathscr A$ has at most $n$ 
members. Thus $\mathcal C_{n+1}(\mathscr A)=\varnothing$ and $d_\sigma(\mathscr A)=n$.

(iii)$\Rightarrow$(ii) Put $\mathscr A_j=u_j\mathscr A$. As above, $\mathscr A=\mathscr
A_1\oplus\cdots\oplus\mathscr A_n$ and, since $\mathscr A_j$ has no non-trivial idempotent, $\mathcal
C_2(\mathscr A_j)=\varnothing$ by \cite[Corollary~3.2]{BK26}.

To prove the uniqueness, suppose that $\mathscr A=\mathscr B_1\oplus\cdots\oplus\mathscr B_m$ is another decomposition as in (ii). 
Let $v_j$ be the unit of $\mathscr B_j$. By the implication (ii)$\Rightarrow$(i) we have $m=n$. As shown above, every idempotent of 
$\mathscr A$ is a sum of some of the idempotents $u_1,\ldots,u_n$. Hence $v_j=\sum_{k\in F_j}u_k$ for a non-empty set 
$F_j\subseteq\{1,\ldots,n\}$. If $F_j$ contained two distinct indices $k$ and $l$, then $u_k$ would be a non-trivial idempotent of 
$v_j\mathscr A=\mathscr B_j$, contradicting the assumption $\mathcal C_2(\mathscr B_j)=\varnothing$. Hence each $F_j$ is a 
singleton, which means that each $v_j$ is one of the $u_k$. Since both families sum to $e$, they differ only by a permutation.
\end{proof}

\section*{Declarations}

\subsection*{Funding}
The first author was supported by the Slovenian Research and Innovation Agency program P2-0268.
The second author was supported by the Slovenian Research and Innovation Agency program P1-0222 and grant J1-50002.

\subsection*{Competing interests}
The authors have no relevant financial or non-financial interests to disclose.

\subsection*{Ethics approval and consent}
Not applicable. This study does not involve human participants, animals, or personal data.

\subsection*{Data availability}
No datasets were generated or analysed during the current study.

\subsection*{Code availability}
Not applicable.

\subsection*{Author contributions}
Both authors contributed equally to the conception of the study, the proofs of the results, and the writing of the manuscript. Both read and approved the final manuscript.

\subsection*{AI Declaration}
During the preparation of this work, the authors used Claude (Anthropic, model Claude Opus 5.5). The AI assistant was used to search 
for and check relevant literature, to propose arguments and examples, and for language and \LaTeX\ editing. All mathematical 
arguments were independently checked and revised by the authors, who take full responsibility for the content and the final version of 
the manuscript.


\end{document}